\documentclass[11pt]{article}
\usepackage[utf8]{inputenc} 
\usepackage[T1]{fontenc}    
\usepackage{url}            
\usepackage{booktabs}       
\usepackage{amsfonts}       
\usepackage{nicefrac}       
\usepackage{microtype}      
\usepackage{natbib}
\usepackage[ruled,vlined]{algorithm2e}
\usepackage{mathtools}
\usepackage{amssymb}
\usepackage{amsthm}
\usepackage{amsmath}
\usepackage{verbatim}  
\usepackage{enumerate}
\usepackage{fancyvrb}
\usepackage{bm}
\usepackage{dsfont}
\usepackage{subcaption}
\usepackage{mwe}
\usepackage{caption}
\usepackage{xcolor}
\usepackage[margin=1.4in]{geometry}
\usepackage{appendix}

\begin{document}

\title{BigSurvSGD: Big Survival Data Analysis via Stochastic Gradient Descent}

\author{Aliasghar Tarkhan and Noah Simon \\
Department of Biostatistics\\
       University of Washington\\
       Seattle, WA 98195-4322, USA\\
atarkhan@uw.edu and nrsimon@uw.edu}


\maketitle
\begin{abstract}
In many biomedical applications, outcome is measured as a ``time-to-event'' (eg. disease progression or death). To assess the connection between features of a patient and this outcome, it is common to assume a proportional hazards model, and fit a proportional hazards regression (or Cox regression). To fit this model, a log-concave objective function known as the ``partial likelihood'' is maximized. For moderate-sized datasets, an efficient Newton-Raphson algorithm that leverages the structure of the objective can be employed. However, in large datasets this approach has two issues: 1) The computational tricks that leverage structure can also lead to computational instability; 2) The objective does not naturally decouple: Thus, if the dataset does not fit in memory, the model can be very computationally expensive to fit. This additionally means that the objective is not directly amenable to stochastic gradient-based optimization methods. To overcome these issues, we propose a simple, new framing of proportional hazards regression: This results in an objective function that is amenable to stochastic gradient descent. We show that this simple modification allows us to efficiently fit survival models with very large datasets. This also facilitates training complex, eg. neural-network-based, models with survival data.

\end{abstract}

\section{Introduction}

It is commonly of interest in biomedical settings to characterize the relationship between characteristics of an individual and their risk of experiencing an event of interest \citep[eg. progression of a disease, recovery, death, etc. see][]{Lee1997}. Outcomes of this type are known as ``time-to-event'' outcomes, and characterizing such relationships is known as Survival analysis \citep{Schober2018}. In such applications, we often only have partial information on some patients due to censoring (e.g., they might leave the study before experiencing the event of interest). 

Cox proportional hazards regression (CoxPH) \citep{Cox1972} is the most common tool for conducting survival analyses. The CoxPH model assumes a particular semi-parametric relationship between the risk at each given time of experiencing an event and the features of a patient (eg. age, sex, treatment assignment, etc). To estimate, parameters in this model CoxPH regression maximizes a log-concave function known as the "partial likelihood". Once estimated, this model can provide a person-specific hazard of experiencing an event (as a function of their features). Such predictions are often used in personalized medicine, eg. in the development of prognostic and predictive biomarkers \citep{Bjorn2008}. To maximize the partial likelihood, it is most common to use an efficient second-order algorithm such as Newton-Raphson \citep{MITTAL2014} for datasets with few features (though potentially many observations). Traditionally, CoxPH has been used on data-sets with relatively few features, though penalized extensions have been developed for high dimensional applications \citep{Simon2011}.

It is increasingly common to have biomedical datasets with a large number of observations, especially with increasing use of electronic medical records \citep{Raghupathi2014}. Although the CoxPH regression model has been widely used for small-to-moderate numbers of observations, current methodologies for fitting the Cox model have issues on datasets with many observations. In particular, in fitting the Cox model, it is common to leverage the sequential structure of the partial likelihood to vastly speed up computation \citep[from $O(n^2)$ to  $O(n)$, see][]{Simon2011}. However, when there are a large number of observations, this can lead to computational instability (which we illustrate in this manuscript).



The second issue is that the partial likelihood does not naturally decouple over individuals or subsets of individuals. Thus, if the dataset does not fit in memory, the model can be very computationally expensive to fit: Standard distributed optimization methods such as those based the alternating direction method of multipliers \citep{Boyd2010} cannot be used. This additionally means that the objective is not directly amenable to stochastic gradient-based optimization methods \citep{Ruder2016}. Unfortunately, to fit more complex neural network-based models, it is most common to use stochastic-gradient-based optimization \citep{NN}. This decoupling issue makes it impossible or at least very impractical to fit neural-network-based models with time-to-event data.

 
In this paper, we propose a novel and simple framework for conducting survival analysis using the CoxPH model. Our framework is built upon an objective that is a modification of the usual partial likelihood function. In particular, this modified objective decouples over subsets of observations, and allows us to employ stochastic-gradient-based algorithms that engage only a subset of our data at each iteration. We show that the parameters estimated by this new objective function are equivalent to the original parameters when the assumptions of the CoxPH model hold (and may actually be more robust in the case of model misspecification). In addition, our new objective function is amenable to optimization via stochastic-gradient based algorithms that are computationally efficient and stable for this objective and can easily scale to datasets that are too large to fit in memory. We further discuss how our new framework can be implemented in both streaming \citep{Gaber2005} and non-streaming algorithms. We discuss extending our framework to use mini-batches and we present some recommendations that we have found important, in practice, for performance. We also extend our methodology to produce asymptotically valid confidence intervals.

We organize this paper as follows. In Section~\ref{sec2}, we review related work in the literature. In Section~\ref{sec3}, we first review the CoxPH regression model, the partial likelihood, and standard optimization tools used to maximize the partial likelihood. We then present our new framework for fitting the CoxPH model and prove the statistical equivalence of the parameters indicated by our optimization problem and those from the standard CoxPH model. We also discuss streaming and non-streaming algorithms for parameter estimation in our framework. In addition, we discuss some recommendations that we have found important in practice for performance. 
In Section~\ref{sec5}, we frame our estimator as the minimizer of a U-statistic-based objective function. In Section~\ref{sec6}, we use this framing, and propose two approaches to construct asymptotically valid confidence intervals for our proposed parameters. In Section~\ref{sec7}, we provide simulation results that compare estimates from our proposed framework to the current gold-standard estimator based on the full partial likelihood. In Section~\ref{sec8}, we conclude our paper and discuss some potential implications of our framework. 

\section{Related work}
\label{sec2}
Although there are many fully-parametric \citep[e.g.,][]{Wei1992} and non-parametric \citep[e.g.,][]{Chapfuwa2018} approaches for survival analysis, we focus on those that are based on the Cox proportional hazards model (as those are most relevant to our proposal). \citet{Toulis2017} presented SGD-based algorithms for a variety of applications including the Cox proportional hazards model. However their algorithm suffers from two issues: It cannot accommodate streaming data, and in fact requires $\sim n^2$ computation (which we are aiming to avoid with our proposal). \citet{Raykar2008} proposed directly maximizing the concordance index \citep{Peter2012}. While this is an interesting predictive target, it moves us away from generative parameters in the Cox model. Additionally, as this objective function is discontinuous, confidence intervals may be difficult to obtain. \citet{Katzman2017, Ching2018} connected neural networks to the log-partial likelihood. However, they engaged with [non-stochastic] gradient descent, which as we will discuss in this paper, is not amenable to very large and/or distributed datasets. The work of \citet{Kvamme2019} is most closely related to ours: As with \citet{Katzman2017} and \citet{Ching2018}, they connect neural networks with the partial likelihood; however they note that a stochastic gradient-like optimization method will be needed to scale to large datasets. As such, they come up with a heuristic for an ``approximate gradient''. They do not justify the heuristic (and in fact, their stochastic gradient is not unbiased, so there is no guarantee that any of the results of SGD-based methods will hold). Our reframing of proportional hazards modeling generalizes and justifies their heuristic --- it both identifies why it should work and proves that it will.
Our contributions in this paper are:
\begin{itemize}
    \item We propose a new optimization framework based on a modified version of the full log-partial likelihood that is amenable to stochastic/online optimization, requires $\sim n$ computation, and is not susceptible to floating-point error. 
    \item We prove that the parameter estimated by our framework is equivalent to the generative parameter in a correctly specified CoxPH model. 
    \item We prove that for a linear CoxPH model, our proposed estimator is statistically rate optimal (MSE converges at a rate of $n^{-1}$). 
    \item We theoretically and empirically show that for a linear CoxPH model, our framework accommodates construction of asymptotically valid confidence intervals.
    \item We empirically show that our framework may be more computationally efficient and stable than the current gold standard implementation for fitting the CoxPH model when the dataset is large.
\end{itemize}

\section{Methodology}
\label{sec3}
\subsection{Cox Proportional Hazards (Cox PH) Model}
The Cox PH model, proposed by \citet{Cox1972}, is a commonly used semi-parametric regression model for evaluating the association between the time until some event of interest and a set of variable(s) measured on a patient. More formally, suppose on each patient we measure $T$ an event time, and $X = (X_1,\ldots, X_p)$ a vector of numeric features. The Cox model engages with the so-called hazard function
\[
h(x,t) = \frac{p(t|x)}{S(t|x)}
\]
where $p(t|x) = \frac{d}{dt}P(T < t|X=x)$, $S(t|x) = P(T > t|X=x)$. The hazard function, $h(x,t)$, can be thought of as the probability density of having an event at time $t$, given that a patient (with covariates $x$) has not had an event up until that time. In particular the Cox model assumes a particular form for the hazard function:
\begin{align}
h(t, \boldsymbol{x};\boldsymbol{\beta}^*) = h_0(t)\exp{(f_{\boldsymbol{\beta}^*}(\boldsymbol{x}))}
\label{e1}
\end{align}
where $f_{\boldsymbol{\beta}^*}$ is a specified function of parameters $\boldsymbol{\beta}^*=(\beta^*_1, \beta^*_2, \dots, \beta^*_k)$ that determines the role played by $\boldsymbol{x}$ in the hazard; and $h_0(t)$ is a baseline hazard function (independent of covariates). Note that $f_{\boldsymbol{\beta}^*}(\boldsymbol{x})$ may be assumed to be of different forms in different applications: For instance, in many scenarios $k$ is taken to be $p$, and the simple linear model $f_{\boldsymbol{\beta}}(\boldsymbol{x})=\boldsymbol{x}^T\boldsymbol{\beta}=  + \beta_1 x_1 + \dots + \beta_p x_p$ is used. This model assumes that the manner in which a patient's covariates modulate their risk of experiencing an event is independent of time. In particular it is encoded entirely in $f_{\boldsymbol{\beta}^*}$. This simplifies estimation and interpretation of the predictive model.

Our aim is to use data to estimate $\boldsymbol{\beta}^*$. In particular we will assume that we have a dataset with $n$ independent observations drawn from model \eqref{e1}: $\mathcal{D}^{(n)}= \{\mathcal{D}_i=(y_i, \boldsymbol{x}^{(i)})|i=1,2, \dots, n\}$. For the moment we assume that there is no censoring (all event times are observed), and no ties (all event times are unique). Estimation is conducted using the log-partial-likelihood:
\begin{align}
pl^{(n)}(\boldsymbol{\beta}|\mathcal{D}^{(n)}) = \operatorname{log}\Big(\prod_{i=1}^{n} \frac{h(\boldsymbol{x}^{(i)}; \boldsymbol{\beta})}{\sum_{j \in \mathcal{R}_i}h(\boldsymbol{x}^{(j)}; \boldsymbol{\beta})}\Big) \nonumber\\
&= \sum_{i=1}^{n}\Big(f_{\boldsymbol{\beta}}(\boldsymbol{x}^{(i)})-\operatorname{log}\Big(\sum_{j \in \mathcal{R}_i}\exp{(f_{\boldsymbol{\beta}}(\boldsymbol{x}^{(j)}))}\Big)\Big)
\label{e2}
\end{align}
where $\mathcal{R}_i = \{j\,|\,t_j \geq t_i\}$ is the ``risk set for patient $i$''. Note that aside from $\mathcal{R}_j$, the expression in \eqref{e2} is independent of the event times. Extending this partial likelihood to deal with censoring, left-truncation \citep{Klein2003} and ties is quite straightforward (see Appendix \ref{secAA} of the Supplementary Materials for more details of this extension), however for ease of exposition we do not include it in this manuscript.

\subsection{\protect{Estimation by maximizing the log-partial-likelihood}}
Using the log-partial-likelihood from equation~\eqref{e2} an estimate of $\boldsymbol{\beta}^*$ can be obtained as
\begin{align}
\hat{\boldsymbol{\beta}}^{(n)} = \underset{\boldsymbol{\beta}}{\operatorname{argmin}}\left\{-pl^{(n)}(\boldsymbol{\beta}|\mathcal{D}^{(n)})\right\}
\label{e3}
\end{align}
When linear $f_{\boldsymbol{\beta}}(\bold{x})=\bold{x}^{\top}\boldsymbol{\beta}= \beta_1 x_1 + \dots + \beta_p x_p$ is used, our objective function in \eqref{e3} is convex in $\boldsymbol{\beta}$, and thus the tools of convex optimization can be applied to find $\hat{\boldsymbol{\beta}}^{(n)}$ (see \textit{Corollary 1} in Appendix \ref{secAB} of the Supplementary Materials for the proof of convexity; included for completeness). In the current gold standard {\tt survival} package in {\tt R} \citep{coxph2019}, Newton-Raphson is used to minimize \eqref{e3} with linear $f$. In the later sections of this manuscript we will refer to this implementation as {\tt coxph()}. 

For linear $f$, one can show that $\left\|\hat{\boldsymbol{\beta}}^{(n)} - \boldsymbol{\beta}^{*}\right\|_2^2 = O_p\left(n^{-1}\right)$ which is rate optimal \citep[as is standard for estimation in parametric models, see][]{Van2000}.  

In current state-of-the-art packages, the structure of the ordered structure of the loss (as well as the gradient, and hessian) are leveraged to improve computational efficiency. In particular, we examine the gradient
\begin{align}
    \nabla_{\beta}\Big\{-pl^{(n)}(\boldsymbol{\beta}|\mathcal{D}^{(n)})\Big\}= -\sum_{i=1}^{n}\Big(\dot{f}_{\boldsymbol{\beta}}(\bold{x}^{(i)})-\frac{\sum_{j \in \mathcal{R}_i}\dot{f}_{\boldsymbol{\beta}}(\bold{x}^{(j)})\exp{(f_{\boldsymbol{\beta}}(\bold{x}^{(j)})})}{\sum_{j \in \mathcal{R}_i}\exp{(f_{\boldsymbol{\beta}}(\bold{x}^{(j)}))}}\Big).
    \label{e4}
\end{align}
where $\dot{f}_{\boldsymbol{\beta}}(\bold{x})=\nabla_{\beta}\{f_{\boldsymbol{\beta}}(\bold{x})\}$ is the gradient of $f_{\boldsymbol{\beta}}(\bold{x})$ with respect to $\boldsymbol{\beta}$. While a naive calculation would have $n^2$ computational complexity because of the nested summations, this is not necessary. In the case that the times are ordered $t_1 < t_2 <\ldots < t_n$ we see that $R_{i} = R_{i+1}\cup \{i\}$. This allows us to use cumulative sums and differences to calculate the entire gradient in $O(n)$ computational complexity, with a single $n\operatorname{log}(n)$ complexity sort required at the beginning of the algorithm \citep{Simon2011}. This is also true for calculating the Hessian. Unfortunately, however, when employing this strategy, the algorithm becomes susceptible to roundoff issues, especially with a larger number of observations ($n$) and features ($p$), as seen in Section~\ref{sec:sim-big}.

Additional inspection of the gradient in \eqref{e4} shows why stochastic-gradient-based methods cannot be used to decouple gradient calculations over observations in our sample: While the gradient can be written as a sum over indices $i=1,\ldots, n$, the denominator for the $i=1$ term involves all observations in the dataset. In the next section, we propose a novel simple modification of optimization problem~\eqref{e3} that admits an efficient stochastic-gradient-based algorithm for estimating $\boldsymbol{\beta}^*$.

\subsection{Estimation using SGD: BigSurvSGD}
We begin by reformulating our problem. We consider a population parameter $\boldsymbol{\beta}^{(s)}$, defined as the population minimizer of the expected partial likelihood of $s$ random patients (which we will refer to as ``strata of size s'')
\begin{align}
\boldsymbol{\beta}^{(s)} = \underset{\boldsymbol{\beta}}{\operatorname{argmin}}\Big\{\mathbb{E}_{s}[-pl^{(s)}(\boldsymbol{\beta}|\mathcal{D}^{(s)})]\Big\}
\label{e5}
\end{align}
Here we think of $\mathcal{D}^{(s)}$ as a draw of $s$ random patients from our population. Note that the minimum value for $s$ is 2, otherwise, expression (2) becomes zero for all $\boldsymbol{\beta}$. By including a superscript $s$ in $\boldsymbol{\beta}^{(s)}$, we note that this parameter may depend on $s$. In fact, when the assumptions of the Cox model hold \eqref{e1} then we have $\boldsymbol{\beta}^{(s)} = \boldsymbol{\beta}^{*}$ for all $s$. The proof of this is quite simple, with details given in Appendix \ref{secAB} of the Supplementary Materials.

To estimate $\boldsymbol{\beta}^{*}$, we select a small fixed $s$ ($s<<n$) and directly apply stochastic gradient descent to the population optimization problem \eqref{e5}. In practice this will amount to calculating stochastic gradients using random strata of size $s$. One may note that for $s$ small, there are on the order of $n^s$ such strata. However, results for stochastic gradient descent indicate that under strong convexity of \eqref{e5}, on the order of only $n$ steps should be required to obtain a rate optimal estimator (converging at a rate of $n^{-1}$ in MSE). See Appendix \ref{secAC} of the Supplementary Materials for the proof of strong convexity of \eqref{e5} and the convergence rate $O(n^{-1})$. 

In Section~\ref{sec:sim-big}, we see that this modification mitigates issues with roundoff error, and allows us to computationally efficiently fit survival models with millions of observations and many features.

\subsubsection{Pairwise concordance ($s=2$)}
An interesting special case is when we choose strata of size $s=2$, and look at pairs of patients. Then, in the case of no censoring, the population minimizer in \eqref{e5}, i.e., $\boldsymbol{\beta}^{(2)}$ maximizes the expectation of the pairwise log-partial likelihood
\begin{align}
    pl^{(2)}(\boldsymbol{\beta}|\mathcal{D}^{(2)}) &= log\Big(\frac{\exp{(f_{\boldsymbol{\beta}}(\bold{x}^{(1)}))}}{\exp{(f_{\boldsymbol{\beta}}(\bold{x}^{(1)}))}+\exp{(f_{\boldsymbol{\beta}}(\bold{x}^{(2)}))}}\Big)1(t_1<t_2)\nonumber\\
    &+log\Big(\frac{\exp{(f_{\boldsymbol{\beta}}(\bold{x}^{(2)}))}}{\exp{(f_{\boldsymbol{\beta}}(\bold{x}^{(1)}))}+\exp{(f_{\boldsymbol{\beta}}(\bold{x}^{(2)}))}}\Big)1(t_2<t_1).
    \label{ep}
\end{align}
This log-partial-likelihood can be thought of as a smoothed version of the standard concordance measure used in the concordance index \citep{Peter2012}. Thus, even when the proportional hazards model does not hold, the parameter $\boldsymbol{\beta}^{(2)}$ maintains a useful interpretation as the population minimizer of the average smoothed concordance index. Also \eqref{ep} is similar to the objective function for conditional logistic regression (CLR) with strata size $s=2$ \citep{Breslow1980}. In the deep learning literature, neural net models with a conditional logistic outcome layers are often referred to as Siamese Neural Networks \citep{Gregory2015}.

\subsubsection{Optimization with SGD}
\label{secSGD}
Suppose that we have $n_s$ independent strata, $D^{(s)}_1, \ldots, D^{(s)}_{n_s}$ each with $s$ independent patients drawn from our population (with $s\geq 2$). For ease of notation, let $I_m$ denote the indices of patients in strata $D^{(s)}_m$ for each $m \leq n_s$.

We first note, that for any $\boldsymbol{\beta}$ we have
\begin{align}
\nabla_{\boldsymbol{\beta}}\mathbb{E}_s\left[pl^{(s)}(\boldsymbol{\beta}|\mathcal{D}^{(s)}_m)\right]= \mathbb{E}_s\left[\nabla_{\boldsymbol{\beta}} \left\{pl^{(s)}(\boldsymbol{\beta}|\mathcal{D}^{(s)}_m)\right\}\right], \qquad \textrm{ for all $m\leq n_s$}
\label{e6}
\end{align}
when $\bold{x}$ are drawn from a reasonable distribution (eg. bounded); and $f_{\boldsymbol{\beta}}(\bold{x})$ is not too poorly behaved (eg. Lipschitz). Here $\nabla_{\boldsymbol{\beta}} \left\{pl^{(s)}(\boldsymbol{\beta}|\mathcal{D}^{(s)}_m)\right\}$ is defined analogously to \eqref{e5} using $\mathcal{D}^{(s)}_m$
\begin{align}
    \nabla_{\beta}\Big\{-pl^{(s)}(\boldsymbol{\beta}|\mathcal{D}^{(s)}_m)\Big\}= -\sum_{i\in I_m} \left(\dot{f}_{\boldsymbol{\beta}}(\bold{x}^{(i)})-\frac{\sum_{j \in \mathcal{R}^m_i}\dot{f}_{\boldsymbol{\beta}}(\bold{x}^{(j)})\exp{(f_{\boldsymbol{\beta}}(\bold{x}^{(j)})})}{\sum_{j \in \mathcal{R}^m_i}\exp{(f_{\boldsymbol{\beta}}(\bold{x}^{(j)}))}}\right)
    \label{e7}
\end{align}
where $R_i^m = \{j\,|\,t_j \geq t_i\, \textrm{ and } i,j\in I_m\}$ are risk sets that include only patients in stratum $m$; and $\dot{f}_{\beta}$ denotes the gradient of $f$ wrt $\beta$.

From here we can give the simplest version of our stochastic gradient descent (SGD) algorithm for \eqref{e3}. We choose an initial $\hat{\boldsymbol{\beta}}(0)$ (perhaps $=0$), and at each iteration $m=1,\ldots, n_s$, we update our estimate by
\begin{align}
\hat{\boldsymbol{\beta}}(m) = \hat{\boldsymbol{\beta}}(m-1) + \gamma_m \times \nabla_{\beta} \left\{pl^{(s)}(\hat{\boldsymbol{\beta}}(m-1)|\mathcal{D}^{(s)}_m)\right\}.
\label{e8}
\end{align} 
Here, $\gamma_m$ is the learning rate (and should be specified in advance, or determined adaptively as discussed later in Section~\ref{str} and Appendix \ref{secAD} of the Supplementary Materials). The computation time to run $n_s$ steps of stochastic gradient descent according to \eqref{e8} for linear $f_{\boldsymbol{\beta}}$ is $\sim sn_s  p = np$ (where $n = sn_s$ is the total sample size). If ordering of risk sets is not leveraged, then $\sim nps$ computation is required. In contrast, Newton's algorithm for optimizing the full log-partial likelihood requires $\sim np^2$ computation per iteration when the roundoff-error-prone updating rule is used (and $\sim n^2p + np^2$ if not). Additionally, using these small strata of size $s$, we are not prone to roundoff issues when using stochastic optimization (because this sum is calculated separately for each strata). 

It has been shown that SGD algorithms for \emph{strongly convex objective-functions} are asymptotically more efficient if we use a running \textbf{average} of the iterates as the final estimate \citep{Ruppert1988, Polyak1992}. Additionally in this case $\gamma_m = \gamma$ can be set to a fixed value (so long as it is sufficiently small). We denote the running-average estimator by 
\begin{align}
\tilde{\boldsymbol{\beta}}(m) = \frac{1}{m}\sum_{i=1}^{m}\hat{\boldsymbol{\beta}}(i)
\label{e9}
\end{align}  
Note that the averaging process does not change the values of $\hat{\boldsymbol{\beta}}(m)$. In our simulations in the main manuscript, we use the averaged $\tilde{\boldsymbol{\beta}}(m)$. Strong convexity of the objective in \eqref{e5} depends on properties of $f_{\boldsymbol{\beta}}$ and (weakly) on the distribution of $\boldsymbol{x}$. For linear $f_{\boldsymbol{\beta}}$, and $\boldsymbol{x}$ with a non-degenerate distribution, this objective will be strongly convex (See Appendix \ref{secAC} of the Supplementary Materials for the proof of strong convexity of our objective function). In such cases, standard results \citep{Bottou2010} show that $\left\|\tilde{\boldsymbol{\beta}}(m) - \boldsymbol{\beta}^{(s)}\right\|_2^2 = O_p(m^{-1})$. As a reminder, this is the statistically optimal rate of convergence for estimating $\boldsymbol{\beta}^{(s)}$ (or equivalently, $\boldsymbol{\beta}^{*}$ when the Cox model holds) from $m$ observations. See Appendix \ref{secAC} of the Supplementary Materials for the proof of this convergence rate $O(m^{-1})$ for averaging over iterates. When we use averaging over iterates, it has been shown that choosing the learning rate as $\gamma_m=\frac{C}{\sqrt{m}}$ (C is a constant) gives us such an optimal convergence rate \citep{Eric2011}. We choose this learning rate in our simulation studies in Section \ref{sec7} where we tune the constant $C$ to get the desired convergence rate.

\subsubsection{Streaming vs non-streaming}
\label{str}
In the discussion above, we imagined that observations were arriving in a continual stream of strata, and that we were more concerned with the cost of computation than the cost of data collection. The algorithm we described engaged with each stratum only once (in calculating a single stochastic gradient). Algorithm 1 in Appendix \ref{secAE1} of the Supplementary Materials details an implementation of a streaming algorithm in this imagined scenario. In practice, we generally have a fixed (though potentially large) number of observations (non-streaming), $n$, that are not naturally partitioned into strata. To employ SGD here, we randomly partition our observations into $n_s$ disjoint strata of size $s$, and then carry out our updates in \eqref{e8}. Additionally, in practice we may want to take more than one pass over the data (and similarly use more than $1$ random partitioning) and use mini-batches to improve the performance. Algorithm 2 in Appendix \ref{secAE2} of the Supplementary Materials details an implementation that takes multiple passes over the data and use mini-batches (we discuss additional bells and whistles in Appendix \ref{secAD} of the Supplementary materials).

In practice, we may also want to modify the learning rate over iterates to improve the performance of SGD algorithm. We found using algorithm {\tt AMSGrad} \citep{Sashank2019} for adaptively selecting the learning rate using moment, using average over iterates as we discussed in Section \ref{secSGD}), and having $\sim 100$ epochs of data substantially improves performance (see Appendix \ref{secAD} of the Supplementary Materials for more details).

\section{Equivalence to U-statistic based optimization}
\label{sec5}
While in this manuscript we discuss obtaining an estimator by directly attempting to minimize the population objective function \eqref{e5} using SGD, there is a corresponding empirical minimization problem. In particular, for strata of size $s$, one could define an estimator $\hat{\beta}^s$ by
\begin{equation}
\hat{\boldsymbol{\beta}}^{(s)} = \underset{\boldsymbol{\beta}}{\operatorname{argmin}}\Big\{-\sum_{\mathcal{D}^{(s)} \subset \mathcal{D}^{(n)}}pl^{(s)}(\boldsymbol{\beta}|\mathcal{D}^{(s)})\Big\}
\label{e10}
\end{equation}
As in a standard U-statistic \citep{Van2000}, this sum is taken over all subsets of $s$ patients out of our original $n$ patients (resulting in ${n\choose s} \sim n^s$ terms). This appears to be a difficult optimization problem, given the enormous number of terms. However, our approach shows that, in fact, only $\sim n$ of those terms need be considered: The majority contain redundant information. In fact, one could see this directly by noting that the objective function in~\eqref{e10} decouples over subsets: An application of stochastic gradient descent here would involve sampling strata with replacement, and (assuming this empirical objective function is strongly convex), with averaging over iterates, would converge to a tolerance of $1/n$ after $n$ steps. This approach, based on \textit{incomplete U-statistics} \citep{Gunnar1976} could be taken more generally for losses defined by $U$-statistics.

\section{Inference}
\label{sec6}
Reframing our estimator as the minimizer of an $s$-order U-statistic, as in~\eqref{e5}, is quite useful for developing confidence intervals. The large sample properties of such estimators have been studied in the literature \citep{Honore1994, Sherman1994}: In particular these estimators are asymptotically linear \citep{Honore1994}, with an influence function based on their Hoeffding projection \citep{Van2000}. There are two standard approaches to calculating a CI using such a U-statistic-based estimator. The first uses a plug-in estimate of the analytically calculated influence function to compute an estimated variance. The second approach is to use a non-parametric bootstrap. In the following subsections, we briefly discuss these two approaches.

\subsection{Plug-in based approach for estimating CI}
\label{sec6.1}
We note, in optimization problem~\eqref{e10}, the kernel $-pl^{(s)}(\boldsymbol{\beta}|\mathcal{D}^{(s)})$ is differentiable for all values of $\boldsymbol{\beta}$, and that $\boldsymbol{\beta}^{(s)}$ is the unique minimizer of $E[-pl^{(s)}(\boldsymbol{\beta}|\mathcal{D}^{(s)})]$. \citep{Honore1994} studied the asymptotic properties of approximate minimizers of problems \eqref{e10} that satisfy [a generalization of] the approximate first-order term condition 
\begin{align}
    {n\choose s}^{-1} \sum_{\mathcal{D}^{(s)} \subset \mathcal{D}^{(n)}}\frac{\partial}{\partial \boldsymbol{\beta}}\{-pl^{(s)}(\hat{\boldsymbol{\beta}}|\mathcal{D}^{(s)})\} = o_p(n^{-1/2}).
    \label{u2}
\end{align}
Under weak regularity conditions, their work shows that the asymptotic behaviour of $\hat{\boldsymbol{\beta}}$ is given by
\begin{align}
    \sqrt{n}(\hat{\boldsymbol{\beta}}-\boldsymbol{\beta}) \rightarrow N(\bold{0}, \bold{H}^{-1}\bold{V}\bold{H}^{-1})
    \label{u3}
\end{align}
where matrices $\bold{V}$ and $\bold{H}$ are given by 
\begin{align}
    \bold{V} &= s^2\operatorname{E}[r(\mathcal{D}, \boldsymbol{\beta})r(\mathcal{D}, \boldsymbol{\beta})^T] \nonumber\\
    \bold{H} &= \frac{\partial}{\partial \boldsymbol{\beta}}\operatorname{E}[r(\mathcal{D}, \boldsymbol{\beta})]
    \label{u4}
\end{align}
with $r(\mathcal{D}, \boldsymbol{\beta})$ in our problem is defined  as
\begin{align}
    r(\mathcal{D}, \boldsymbol{\beta})=\operatorname{E}\left[\frac{\partial}{\partial \boldsymbol{\beta}}\{-pl^{(s)}(\boldsymbol{\beta}|\mathcal{D}^{(s)}=\{\mathcal{D}, \mathcal{D}_{i_2}, \dots, \mathcal{D}_{i_s}\})\}\middle|\mathcal{D}\right],
    \label{uD}
\end{align}
the expectation of the first derivative of the kernel with respect to $\boldsymbol{\beta}$ when one observation (i.e., $\mathcal{D}$) out of $s$ observations is kept fixed \citep[for more details, see][]{Honore1994}. Finally, \citep{Honore1994} proposed consistent estimators for matrices $\bold{V}$ and $\bold{H}$ as
\begin{align}
    \hat{\bold{V}}_n &= \frac{s^2}{n} \sum_{i=1}^{n}\hat{r}(\mathcal{D}_i, \boldsymbol{\beta})\hat{r}(\mathcal{D}_i, \boldsymbol{\beta})^T \nonumber\\
    \hat{\bold{H}}_n &= \frac{\partial}{\partial \boldsymbol{\beta}} \left(\frac{1}{n} \sum_{i=1}^{n} \hat{r}(\mathcal{D}_i, \boldsymbol{\beta})\right)
    \label{u5}
\end{align}
with $\hat{r}(\mathcal{D}_i, \boldsymbol{\beta})$ given by
\begin{align}
    \hat{r}(\mathcal{D}_i, \boldsymbol{\beta}) = {n-1\choose s-1}^{-1} \sum_{\mathcal{D}_{i_2}, \dots, \mathcal{D}_{i_s} \subset \{\mathcal{D}^{(n)}\backslash {\mathcal{D}_i} \}}\frac{\partial}{\partial \boldsymbol{\beta}}\{-pl^{(s)}(\boldsymbol{\beta}|\mathcal{D}^{(s)}=\{\mathcal{D}_i, \mathcal{D}_{i_2}, \dots, \mathcal{D}_{i_s}\})\},
    \label{u6}
\end{align}
an empirical estimate of \eqref{uD} where we fix $\mathcal{D}=\mathcal{D}_i$ and we include all possible ${n-1\choose s-1}$ combinations of $s-1$ observations $\{i_2, i_3, \dots, i_s\}$ from $\{1, 2, \dots, n\}\backslash \{i\}$.

This plug-in estimator provides an explicit form for our standard error, however it is prohibitively computationally expensive because we need to go through ${n-1\choose s-1}$ distinct combinations for each of our $n$ observations. Our solution is to only consider a random small fraction of ${n-1\choose s-1}$ combinations (e.g., $n_o=1000 << {n-1\choose s-1}$ distinct combinations) for each observation. Algorithm 3 in Appendix \ref{secAF1} of the Supplementary Materials summarizes our implementation of this plugin approach to estimate a $(1-\alpha)\times100\%$ CI. In addition in Section~G of the Supplementary Material, we give an empirical evaluation of the number of sampled strata required for near nominal coverage.

\subsection{Bootstrap based approach for estimating CI}
\label{sec6.2}
The Bootstrap \citep{Efron1979, Bickel1981} is a popular approach for constructing asymptotically valid confidence intervals that requires no analytical work. For i.i.d. observations $\mathcal{D}^{(n)}=\{\mathcal{D}_1, \mathcal{D}_2, \dots, \mathcal{D}_n\}$, in simpler examples \citep{Bickel1981} suggest using a bootstrap sample (sampling with replacement) $\mathcal{D}^{*(n)}=\{\mathcal{D}_1^*, \mathcal{D}_2^*, \dots, \mathcal{D}_n^*\}$ and calculating the bootstrapped U-statistic as
\begin{align}
    U_n^*(\bold{\beta}) = {n\choose s}^{-1} \sum_{\mathcal{D}^{(s)} \subset \mathcal{D}^{*,(n)}}\{-pl^{(s)}(\boldsymbol{\beta}|\mathcal{D}^{(s)})\}.
    \label{u7}
\end{align}
Let $\hat{\boldsymbol{\beta}}$ denote the estimate given by minimizing the original U-statistic in \eqref{e10} (using the original sample) and $\hat{\boldsymbol{\beta}}^b$ denote the estimate given by minimizing \eqref{u7} with the $b$-th bootstrap sample . We define $\tilde{\boldsymbol{\beta}}_n=\frac{1}{n}\sum_{i=1}^{n}\hat{\boldsymbol{\beta}}$ and $\tilde{\boldsymbol{\beta}}^b_n=\frac{1}{n}\sum_{i=1}^{n}\hat{\boldsymbol{\beta}}^b$ as the Ruppert-Polyak averaged estimates of our original sample and the bootstrap resamples, respectively. \citet{Fang2018} showed that  $\sqrt{n}(\tilde{\boldsymbol{\beta}} - \boldsymbol{\beta}^{*})$ and [the random measure] $\sqrt{n}(\tilde{\boldsymbol{\beta}}^b - \tilde{\boldsymbol{\beta}}) |\mathcal{D}^{(n)}$ converge to the same limiting distribution. Thus, we can construct a confidence interval for $\tilde{\boldsymbol{\beta}}$ using quantiles of $\sqrt{n}(\tilde{\boldsymbol{\beta}}^b - \tilde{\boldsymbol{\beta}})|\mathcal{D}^{(n)}$. In our implementation, we use the estimate $\hat{\boldsymbol{\beta}}$ from the original sample as the initial iterate of the SGD algorithm to calculate the coefficients $\hat{\boldsymbol{\beta}}^b$ for each bootstrap resample. Algorithm 4 in Appendix \ref{secAF2} of the Supplementary Materials summarizes the implementation of a bootstrap approach to calculate a $(1-\alpha)100\%$ CI.

\section{Simulation Study}
\label{sec7}
Here we run simulation studies to empirically evaluate the behaviour of our estimator.
\subsection{Data Generation}
\label{DG}
We generate an event time that follows the Cox PH model \eqref{e1} with simple linear $f_{\boldsymbol{\beta}}$ detailed below. We generate the baseline hazard $h_0(t)$ using an exponential distribution with parameters $\lambda=1$. We generate the censoring and event times independently. The details of the data simulation procedure are given below \citep{Bender2005} 
\begin{align*}
&\boldsymbol{\beta^*}=\boldsymbol{1}_{P\times 1}\\
&X_i \sim Uniform(-\sqrt{3},\sqrt{3})\,\, (\text{unit-variance uniform variable}),\\
&y_i \sim exp(\mu=\exp{(-\boldsymbol{X}_i^T\boldsymbol{\beta^{*}}}))\,\, (\text{time to event/censoring}),\\
&\delta_i \sim Bernoulli(p=1-p_c),\,\, p_c = Pr(t_i>c_i)
\label{e12}
\end{align*}
where $y_i=min(t_i, c_i)$, i.e., time to event or censoring whichever comes first. Here $p_c$, the probability of censoring, is a parameter we can tune. Though this is not written in the form of \eqref{e1}, it is still consistent with the Cox PH model assumptions, with $f_{\boldsymbol{\beta^{*}}}(\boldsymbol{x}) = {\boldsymbol{\beta^{*}}}^{\top}\boldsymbol{x}$. In all comparisons, we include the performance of {\tt coxph()} the gold standard {\tt R} implementation of Newton's algorithm for maximizing the partial likelihood.

\subsection{Small Data Results}
\begin{figure}[p]
    \centering
    \begin{subfigure}[b]{0.496\textwidth}
            \centering
            \includegraphics[width=\textwidth]{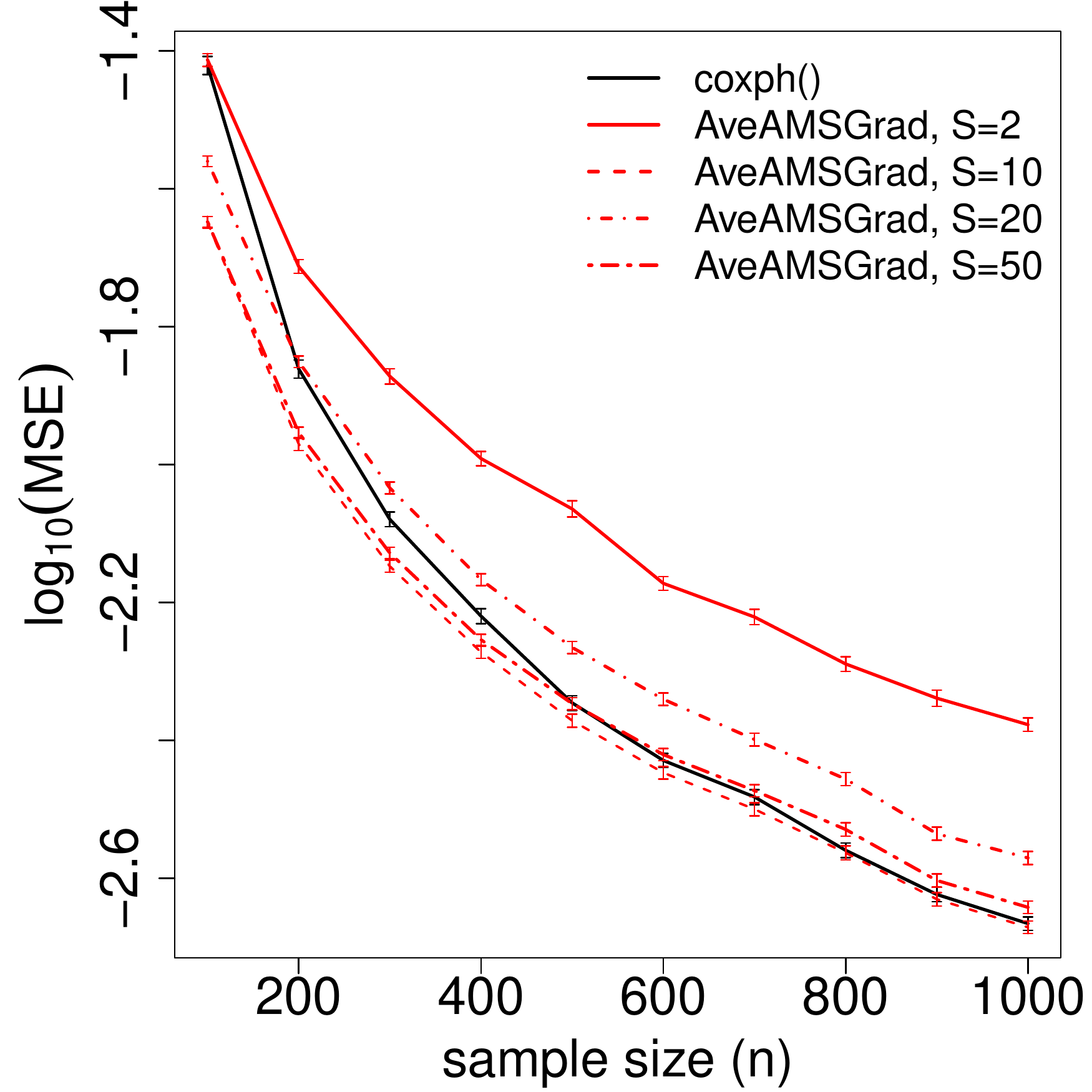}
    \end{subfigure}
    \begin{subfigure}[b]{0.496\textwidth}  
        \centering 
        \includegraphics[width=\textwidth]{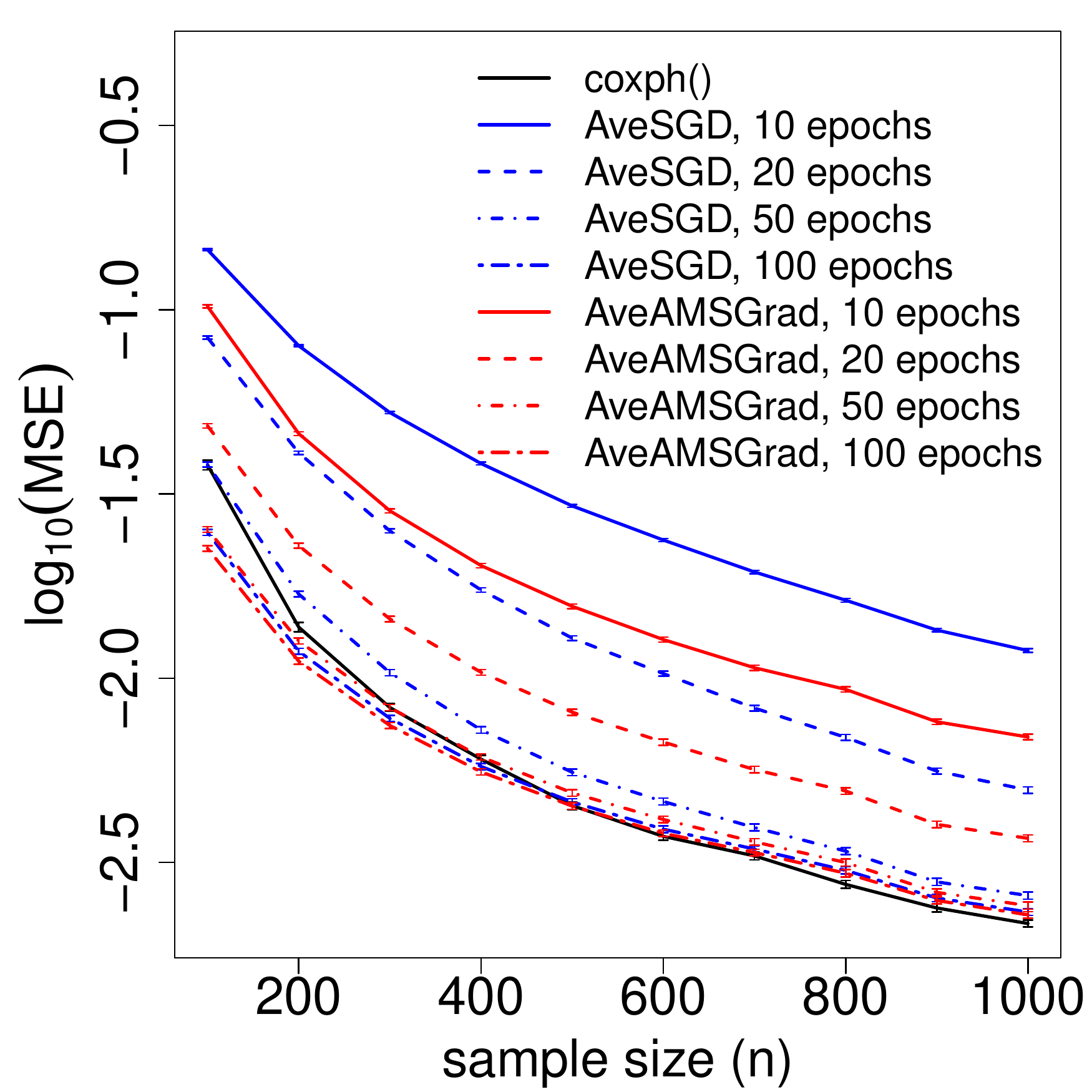}
    \end{subfigure}
    \begin{subfigure}[b]{0.496\textwidth}   
        \centering 
        \includegraphics[width=\textwidth]{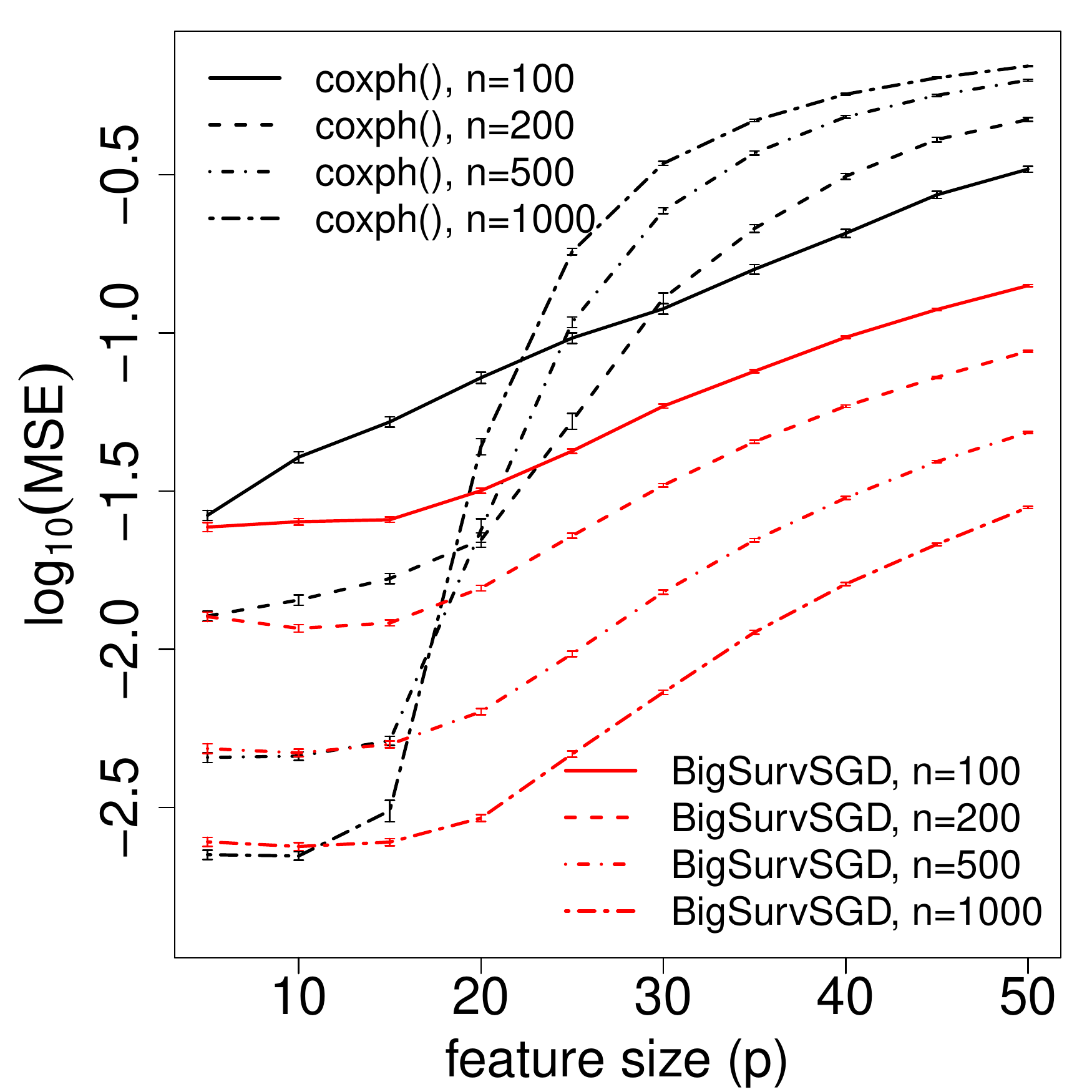}
    \end{subfigure}
    \begin{subfigure}[b]{0.496\textwidth}   
        \centering 
        \includegraphics[width=\textwidth]{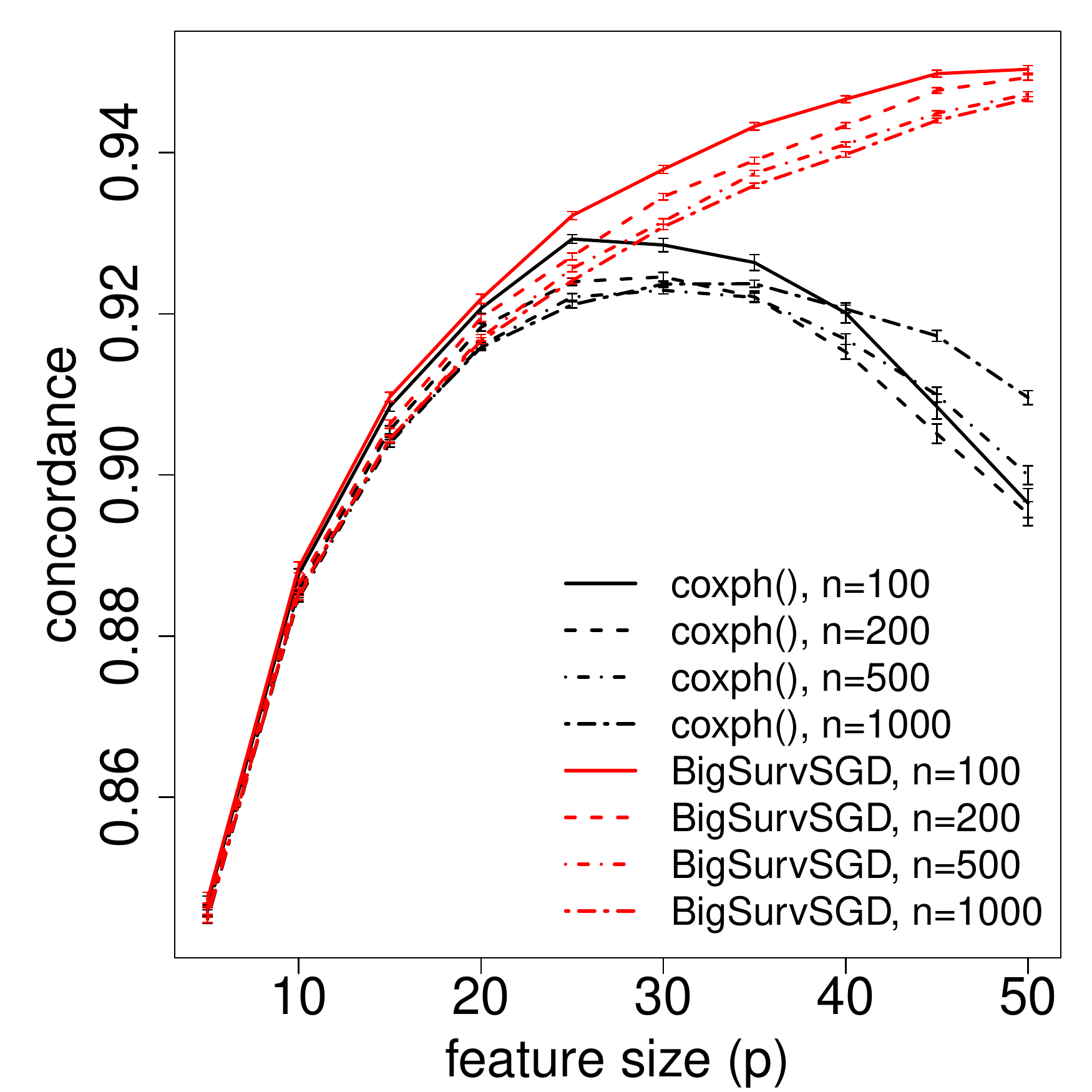}
    \end{subfigure}
    \caption[ The average and standard deviation of critical parameters ]
    {\small (top left) $\operatorname{\log}_{10}$(MSE) of estimates with {\tt{coxph()}} and AveAMSGrad for $P=10$ features, 100 epochs, and different strata sizes ($S$); (top right) $\operatorname{\log}_{10}(MSE)$ of estimates with averaged SGD (AveSGD), averaged AMSGrad (AveAMSGrad), and {\tt{coxph()}} for strata size $S=20$ and $P=10$ features; (bottom left) $\operatorname{log}_{10}$(MSE) and (bottom right) concordance index with BigSurvSGD (AveAMSGrad) and {\tt{coxph()}} for strata size $S=20$, different feature sizes $P$, and different sample sizes $n$. For all results, we choose mini-batch size $K=1$, probability of censoring $p_c=0.2$, and the optimal value $C$ (0.12 for AveAMSGrad and 1.5 for AveSGD) for the learning rate defined as $\gamma_m=\frac{C}{\sqrt{m}}$.} 
    \label{fig01}
\end{figure}

We first evaluate the statistical efficiency of our estimation procedure (using strata sizes of less than $n$). We evaluate mean-squared error (MSE) between estimators $\tilde{\beta}$ (i.e., average over iterates), and $\hat{\beta}^{(n)}$ (i.e., coxph()) and the truth, $\beta^{*}$ over 1000 simulated datasets with up to 100 epochs run. The top left panel in Figure~\ref{fig01} illustrates MSE of $\tilde{\beta}$ with an adaptive learning rate (using {\tt AMSGrad}) for different strata sizes $S$. Although the convergence rate is still $\frac{1}{n}$ for all strata sizes, we see that for small strata sizes (e.g., $S=2$), there is some statistical inefficiency. However, there is nearly no statistical inefficiency for the larger strata sizes (strata sizes larger than $S=20$).

We next evaluate the performance of averaged SGD with a fixed learning rate, against averaged SGD with an adaptive learning rate (using {\tt AMSGrad}) with a fixed strata-size of $S=20$ for both. In addition, we try various numbers of epochs (from 10 to 100). The top right panel in Figure~\ref{fig01} shows performance over 1000 simulated datasets. We see that with enough epochs (around 100) both perform well. However, {\tt AMSGrad} nearly reaches that performance with as few as 50 epochs, where using a fixed learning rate does not attain that performance with fewer than $100$ epochs. For both of these methods, we tuned our [initial] learning-rate to be empirically optimal in these experiments. Hereafter, we use averaged SGD with {\tt AMSGrad} and we simply refer it as {\tt BigSurvSGD}.

In Section \ref{secSGD}, we discussed the computational instability of {\tt coxph()} in small-to-moderate sized datasets and how our framework can avoid such instability. Here we empirically verify those claims. The bottom left and bottom right panels in Figure~\ref{fig01} compare the MSE and concordance index of {\tt{coxph()}} and AveAMSGrad algorithms for small-to-moderate sample sizes ($n$) and varying number of features ($P$) over 1000 simulated datasets. As we see, {\tt{coxph()}} performs poorly for larger $P$ and $n$. For instance, as we see in the bottom left panel in Figure~\ref{fig01}, {\tt{coxph()}} with $(P=50, n=1000)$ performs worse than $(P=50, n=100)$. This is because of computational instability with {\tt{coxph()}} for larger sample sizes and numbers of features. One important aspect of these examples is that we include a large amount of signal (which increases as the number of features increases). With less signal, this instability is less pronounced unless very large sample sizes are used.

\subsection{Big Data Results} \label{sec:sim-big}
\begin{figure}[!tb]
    \centering
    \begin{subfigure}[b]{0.496\textwidth}
            \centering
            \includegraphics[width=\textwidth]{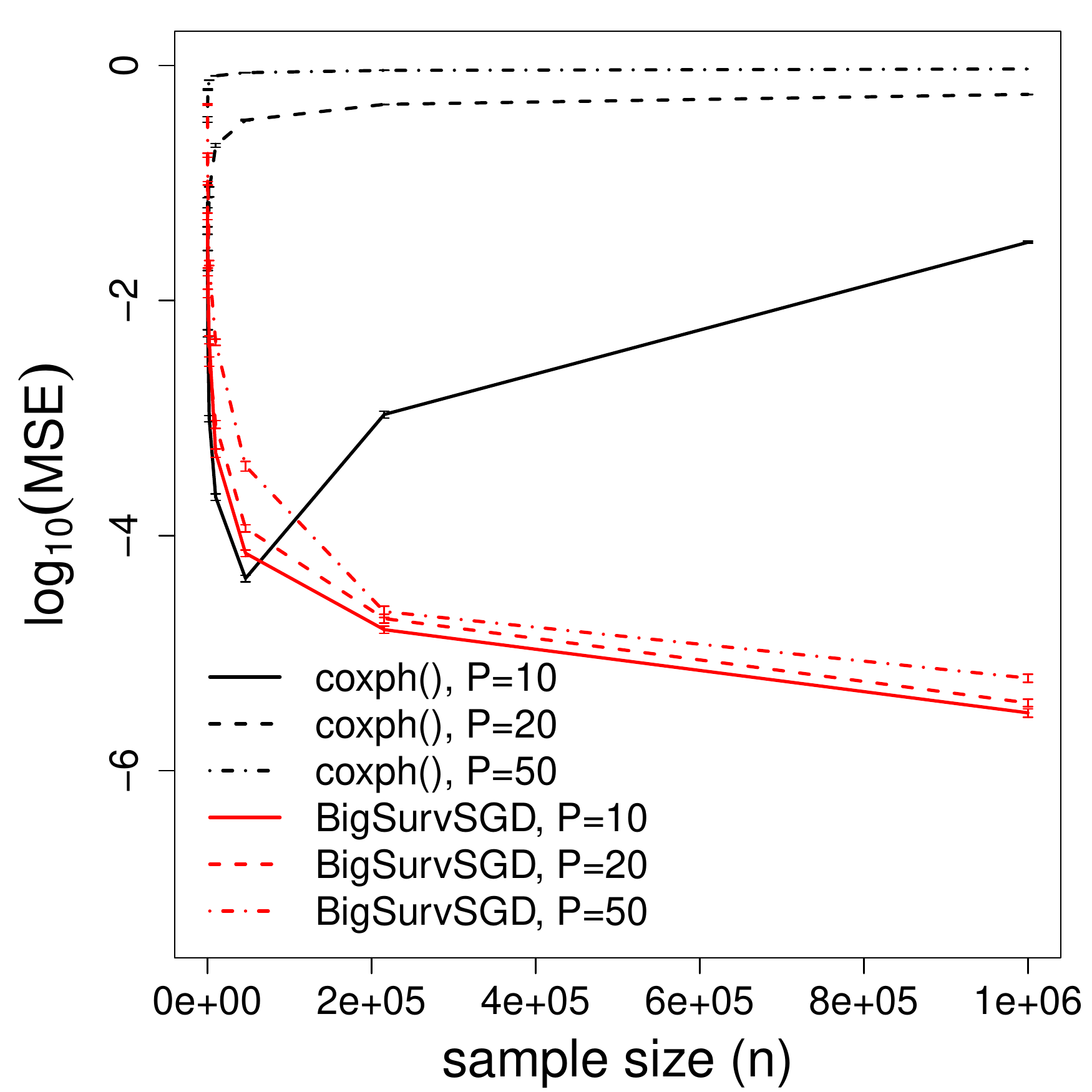}
    \end{subfigure}
    \begin{subfigure}[b]{0.496\textwidth}  
        \centering 
        \includegraphics[width=\textwidth]{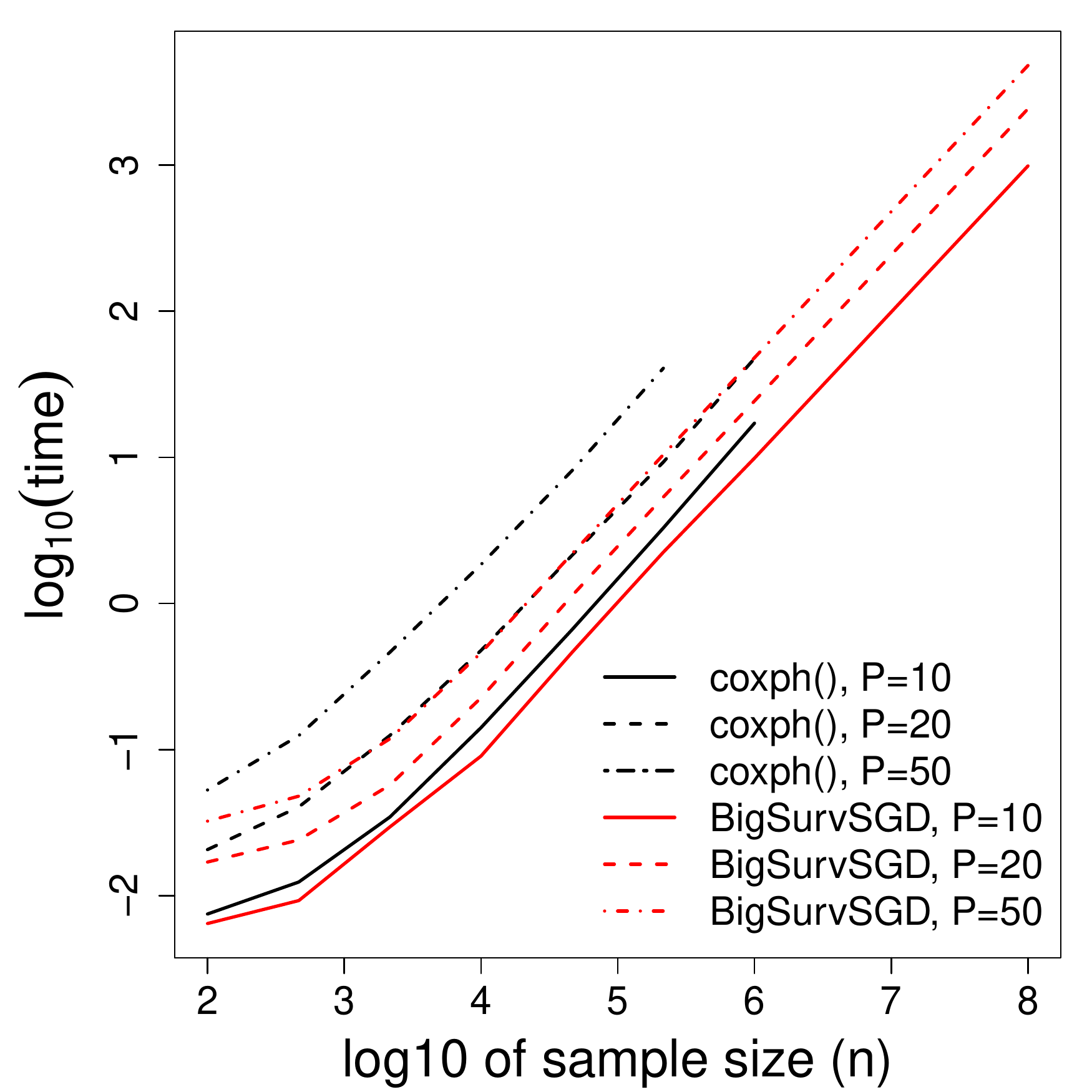}
    \end{subfigure}
    \caption[ The average and standard deviation of critical parameters ]
    {\small (left panel) $\operatorname{\log}_{10}$(MSE) and (right panel) computing time with BigSurvSD (AveAMSGrad) and {\tt{coxph()}} for different sample sizes ($n$) and feature sizes ($P$). We choose mini-batch size $K=1$, probability of censoring $p_c=0.2$, and the optimal value $C=0.12$ for the learning rate defined as $\gamma_m=\frac{C}{\sqrt{m}}$.} 
    \label{fig02}
\end{figure}

We next consider the numerical stability of our algorithm/framework (versus directly maximizing the full partial likelihood using Newton's algorithm). We generated 100 datasets and we only used three epochs for the AveAMSGrad algorithm. Th left panel in Figure~\ref{fig02} shows a surprising and unfortunate result for {\tt coxph()}: We see that as sample size increases drastically, the performance of {\tt coxph()} starts getting worse! In particular, for $P = 20$, {\tt coxph()} is basically producing nonsense by the time we get to $10,000$ observations for this simulation setup. This indicates that for large datasets the current gold standard may be inadequate, though we do note that there is a large amount of signal in these simulations (more than we might often see in practice). In contrast, {\tt BigSurvSGD} has no such issues and gives quite strong performance.

We next examine the computational efficiency of our framework and algorithm. The right panel in Figure~\ref{fig02} plots computing time (in seconds) with {\tt{coxph()}} and BigSurvSGD (AveAMSGrad) for different sample sizes (we did not add error bars for ease of illustration). We considered three different numbers of covariates $P=10, 20$ and $P=50$ to examine the sensitivity of the computing time to the dimension of $\boldsymbol{\beta}$. In these examples, our algorithm read the data in chunks from the hard-drive (allowing us to engage with datasets difficult to fit in memory). The computing time of our proposed framework increases linearly in the sample size, $n$, (and in $p$). The computing time for {\tt{coxph()}} also grows roughly linearly in $n$, though it grows quadratically in $p$ (which can be seen from the poor performance with $p= 50$). Furthermore, {\tt{coxph()}} fails for the medium-to-large datasets as it is poorly equipped to deal with datasets that do not easily fit in memory ({\tt R} unfortunately generally deals somewhat poorly with memory management). For instance, from the right panel in Figure~\ref{fig02}, {\tt{coxph()}} with $P=50$ {\tt coxph()} crashes after $10^5$ observations. As a reminder, the statistical performance of the output of {\tt coxph()}, due to floating-point issues, degenerates much earlier.

\subsection{Inference} \label{sec:inference}
\begin{figure}[!tb]
    \centering
    \begin{subfigure}[b]{0.496\textwidth}
            \centering
            \includegraphics[width=\textwidth]{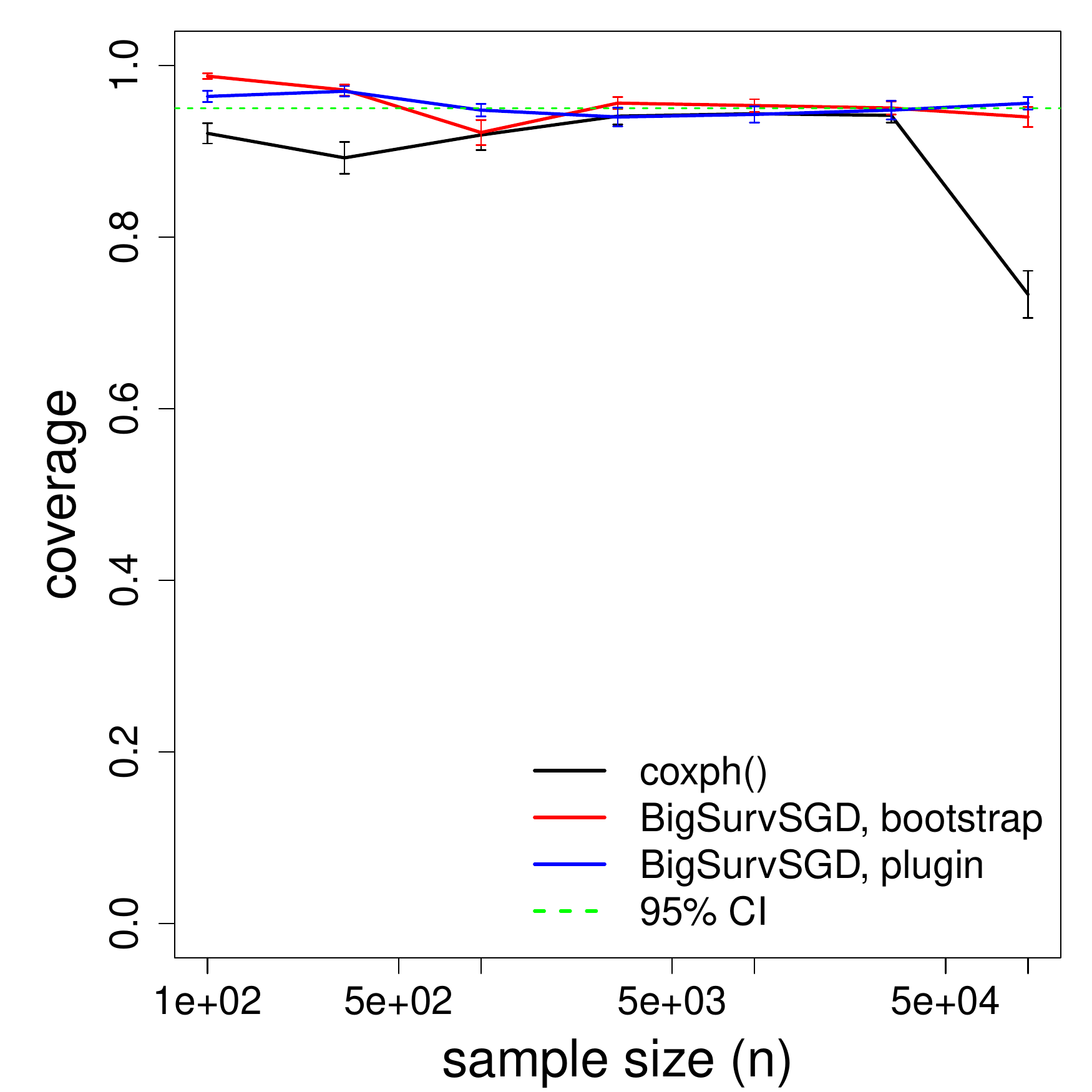}
    \end{subfigure}
    \begin{subfigure}[b]{0.496\textwidth}  
        \centering 
        \includegraphics[width=\textwidth]{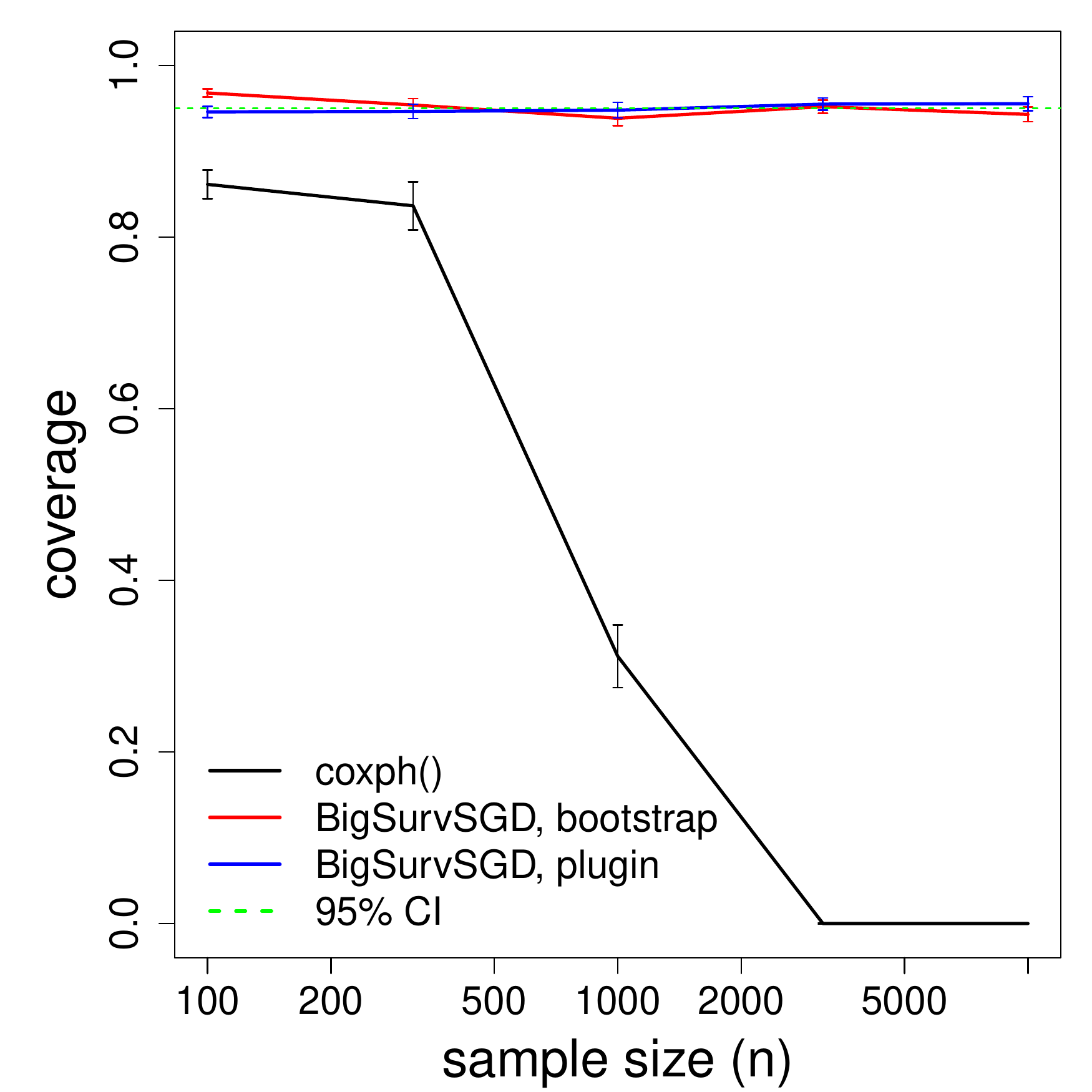}
    \end{subfigure}
    \caption[ The average and standard deviation of critical parameters ]
    {\small (left) coverage for $P=10$ features; (right) coverage for $P=20$ features. For both results, we choose mini-batch size $K=1$, strata size $S=20$, probability of censoring $p_c=0.2$, and the optimal $C=0.12$ for the learning rate defined as $\gamma_m=\frac{C}{\sqrt{m}}$. We use $B=1000$ bootstrap resamples for the bootstrap approach and $n_o=1000$ sample strata per observation for the plugin approach.} 
    \label{fig03}
\end{figure}

We next examine how the coverage of the $95\%$ CIs generated from {\tt{coxph()}}, and our bootstrap and plugin approaches behave as a function of sample size ($n$) and feature size ($P$). We choose strata size $S=20$ and we allow our algorithm to run for 100 epochs to estimate $\tilde{\boldsymbol{\beta}}$. For the bootstrap approach, we consider $B=1000$ bootstrap resamples and we allow each resample to run for 100 epochs to get the bootstrapped estimates $\tilde{\boldsymbol{\beta}}^b$, $b=1, 2, \dots, B=1000$. For the plugin method, we choose $n_o=1000$ sample strata per observation to calculate the standard error of our estimate $\tilde{\boldsymbol{\beta}}$. The left panel in Figure~\ref{fig03} presents the coverage from different approaches for $P=10$ features and varying sample sizes over 100 simulated datasets. We observe all methods perform well for the small-to-medium sample sizes but {\tt{coxph()}} fails for the large sample sizes (e.g., $n=10^5$) due to its computational instability. The right panel in Figure \ref{fig03} presents the coverage from different approaches for more features $P=20$ (higher signal) and varying sample sizes over 100 simulated datasets. We observe both our proposed plug-in and bootstrap approaches give coverage very close to the nominal $95\%$. {\tt{coxph()}} fails for all ranges of sample size; it behaves worse for larger sample sizes.\\

We include additional results from additional simulation experiments as well as data analyses in Appendix \ref{secAG} of the Supplementary Materials. The results are quite similar to those found here.

\section{Conclusion}
\label{sec8}
We propose a simple and novel framework for conducting survival analysis using a Cox proportional hazards model. Our framework leverages a modified optimization problem which allows us to apply iterative methods over only a subset of our observations at a time. In particular, it allows us to leverage the tools of stochastic gradient descent (and its extensions). This results in an algorithm that is more computationally efficient and stable than the current state of the art for problems with a larger number of observations (and features). We also introduced methods to construct confidence intervals for the parameters estimated by our framework. We showed that our framework can handle large survival datasets with little difficulty. This framework will also facilitates training complex models such as neural networks with large survival datasets.

\section*{Acknowledgments}

The authors were funded by NIH grant DP5OD019820 to complete this work.\\
{\it Conflict of Interest}: None declared.

\begin{description}

\item[R-package for  BigSurvSGD:] Our implementation of the SGD procedure is publicly available in the github repository: {\tt https://github.com/atarkhan/bigSurvSGD} \citep{Tarkhan2020}.
\item[Datasets:] We use both simulated experiments (explained in Section \ref{DG}) and datasets FLCHAIN, VETERAN, and GBSG (see Appendix \ref{secAG} of Supplementary Materials for sources of these datasets) (.csv file)

\end{description}

\section*{Appendix}
\appendix
\numberwithin{equation}{section}

\section{Extensions for Left Truncation and Right Censoring}
\label{secAA}
In practice it is common for participants in a study to leave the study before experiencing an event. This is known as right censoring. In particular it is often assumed that a patient has some random ``time until censoring'' $C$, and what we observe is $\operatorname{min}(C,T)$, whichever happens first (the event or censoring), along with $\delta = I(T<C)$ an indicator that the patient experienced an event (rather than censoring). Censored patients still contribute some information to estimation of $\boldsymbol{\beta}^{*}$ --- in particular, if a patient is censored quite late, then there is a long period of time during which we know they did not experience an event (so likely we should estimate them to be low risk). In some studies it is common to consider the event time $T$ as the age at which a patient had an event (rather than the calendar on study). This means that patients do not enroll in the study at $T=0$ (and patients can only be observed to fail once they are enrolled). This phenomenon, wherein patients enroll in a study at times other than $T=0$, is known as left-truncation.

When i) the assumptions of the Cox model hold; and ii) censoring and truncation times are independent of event times conditional on covariates $\boldsymbol{x}$, it is relatively straightforward to adapt cox regression to accommodate these missingness mechanisms. The partial likelihood is modified in 2 minor ways: 1) The outer summation is only taken over indices for which an event occurred (i.e. that were not censored); and 2) The risk sets, $R_i$, are modified to include only patients currently at risk (who have already been enrolled, and have not yet been censored, or had an event) \citep{Andersen1993}. We can similarly modify our objective function \eqref{e5} in Section 3, and apply our algorithm with only minor modification (a slight change in the gradient).

\section{Proof of consistency for parameter $\mathbf{\beta}^{(s)}$ in section 3}
\label{secAB}



In this appendix, for the sake of completeness, we prove the Fisher consistency of parameter $\boldsymbol{\beta}^{(s)}$. We treat the Cox proportional hazard model as a counting process \citep{SurvPoint2008}. We assume that censoring and survival times are independent given the covariate vector of interest $\boldsymbol{x}$ and they follow model \eqref{e1} with true parameter $\boldsymbol{\beta}^*$. We first define some terminology before proceeding with the proof. 

\textit{\textbf{Definition 1: $dN_i(u)$}}. For patient $i$ with time to event $t_i$ define the counting process $dN_i(u)$ by
\begin{align*}
    \int_a^b g(u)dN_i(u) = \Bigg\{ {0\quad\quad\,\,\, if\quad t_i\not\in [a,b] \atop g(t_i)\quad if\quad t_i\in [a,b]}.
\end{align*}
For instance, if we define $g(u)=1$, the above expression is an indicator representing whether patient $i$ failed in interval $[a,b]$ (i.e., 1 represents failure and 0 otherwise). We further define $dN^{(s)}(u)=\sum_{i=1}^s dN_i(u)$ which is a counting process for failure times over all $s$ patients. We assume that the failure time process is absolutely continuous w.r.t. Lebesgue measure on time so that there is at most one failure at any time $u$ (i.e., no ties).

\textit{\textbf{Definition 2: $M_i(u)$}}. We define $M_i(u)$ to be an indicator representing whether patient $i$ is at risk at time $u$, i.e., $t_i \geq u$. By this definition, $M(u)=\sum_{i=1}^sM_i(u)$ indicates number of patients who are at risk at time $u$. Note that the independent censoring assumption implies that those $M(u)$ patients at risk at time $u$ (who have not yet failed or been censored) represent a random sample of the sub-spopulation of patients who will survive until time $u$.

\textit{\textbf{Definition 3: $F(u)$}}. Let $F(u)$ denote the filtration that includes all information up to time $u$, i.e.,
\begin{align*}
    F(u)=\{(dN_i(t), M_i(t), x^{(i)}), i=1, \dots, s\,\quad \text{for $t<u$ and $dN^{(s)}(u)$}\}
\end{align*}
Note that given $F(u)$, we know whether patient $i$ failed or was censored (i.e., $\delta_i$), when they failed/were censored (i.e., $y_i$), and their covariate vectors $\boldsymbol{x}^{(i)}$. \\

Using the above definitions, one can write the log-partial-likelihood as
\begin{align}
    pl^{(s)}(\boldsymbol{\beta}|\mathcal{D}^{(s)}) = \sum_{i=1}^s \int_0^\tau \Big\{f_{\boldsymbol{\beta}}(\boldsymbol{x}^{(i)})-log\Big(\sum_{l=1}^s M_l(u)exp(f_{\boldsymbol{\beta}}(\boldsymbol{x}^{(l)}))\Big)\Big\} dN_i(u)
\end{align}
where $\tau$ is the duration of the study. Note that we keep $f_{\boldsymbol{\beta}}(\boldsymbol{x})$ in the most general form and we only assume that it is differentiable in $\boldsymbol{\beta}$ for all $\boldsymbol{x}$. Then the score function may be written as
\begin{align}
    U^{s}(\boldsymbol{\beta}) &= \frac{\partial pl^{(s)}(\boldsymbol{\beta}|\mathcal{D}^{(s)})}{\partial \boldsymbol{\beta}} \nonumber\\
    &= \sum_{i=1}^s \int_0^\tau \Big(\dot{f}_{\boldsymbol{\beta}}(\boldsymbol{x}^{(i)})-\frac{\sum_{l=1}^s M_l(u)\dot{f}_{\boldsymbol{\beta}}(\boldsymbol{x}^{(l)})exp(f_{\boldsymbol{\beta}}(\boldsymbol{x}^{(l)}))}{\sum_{l=1}^s M_l(u)exp(f_{\boldsymbol{\beta}}(\boldsymbol{x}^{(l)}))} \Big)dN_i(u)\nonumber\\
    &= \sum_{i=1}^s \int_0^\tau \Big(\dot{f}_{\boldsymbol{\beta}}(\boldsymbol{x}^{(i)})-\sum_{l=1}^{s}w_l\dot{f}_{\boldsymbol{\beta}}(\boldsymbol{x}^{(l)}) \Big)dN_i(u)\nonumber\\
    &= \sum_{i=1}^s\int_0^\tau \dot{f}_{\boldsymbol{\beta}}(\boldsymbol{x}^{(i)})dN_i(u)-\sum_{l=1}^{s}\int_0^\tau w_l\dot{f}_{\boldsymbol{\beta}}(\boldsymbol{x}^{(l)}) dN^{(s)}(u).
    \label{A2}
\end{align}
where $w_l = \frac{M_l(u)exp(f_{\boldsymbol{\beta}}(\boldsymbol{x}^{(l)}))}{\sum_{l=1}^s M_l(u)exp(f_{\boldsymbol{\beta}}(\boldsymbol{x}^{(l)}))}$ is a weight proportional to the hazard of failure of patient $l$; $dN^{(s)}(u)=\sum_{i=1}^sdN_i(u)$.\\

Now we show that the parameter $\boldsymbol{\beta^{(s)}}$ is Fisher consistent, i.e., $E[U^{s}(\boldsymbol{\beta^*})]=0$. 
\begin{align}
    E_{\boldsymbol{\beta^*}}[U^{s}(\boldsymbol{\beta^*})] &= E_{\boldsymbol{\beta^*}}\Big[\sum_{i=1}^s\int_0^\tau \dot{f}_{\boldsymbol{\beta}}(\boldsymbol{x}^{(i)})dN_i(u)-\sum_{l=1}^{s}\int_0^\tau w_l\dot{f}_{\boldsymbol{\beta}}(\boldsymbol{x}^{(l)})dN^{(s)}(u)\Big]\nonumber\\
    &=\sum_{i=1}^s\int_0^\tau E_{\boldsymbol{\beta^*}}\Big[\dot{f}_{\boldsymbol{\beta}}(\boldsymbol{x}^{(i)})dN_i(u)\Big]-\sum_{l=1}^{s}\int_0^\tau E_{\boldsymbol{\beta^*}}\Big[w_l\dot{f}_{\boldsymbol{\beta}}(\boldsymbol{x}^{(l)})dN^{(s)}(u)\Big]\nonumber\\
    &\stackrel{(a)}{=}\sum_{i=1}^s\int_0^\tau E_{F(u)}\Big[E_{\boldsymbol{\beta^*}|F(u)}[\dot{f}_{\boldsymbol{\beta}}(\boldsymbol{x}^{(i)})dN_i(u)]\Big]\nonumber \\
    &\,\,- \sum_{l=1}^{s}\int_0^\tau E_{F(u)}\Big[E_{\boldsymbol{\beta^*}|F(u)}[w_l\dot{f}_{\boldsymbol{\beta}}(\boldsymbol{x}^{(l)})]dN^{(s)}(u)\Big]\nonumber\\
    &\stackrel{(b)}{=}\sum_{i=1}^s\int_0^\tau E_{F(u)}\Big[E_{\boldsymbol{\beta^*}|F(u)}[\dot{f}_{\boldsymbol{\beta}}(\boldsymbol{x}^{(i)})dN_i(u)]\Big]\nonumber\\ 
    &- \sum_{l=1}^{s}\int_0^\tau E_{F(u)}\Big[w_l\dot{f}_{\boldsymbol{\beta}}(\boldsymbol{x}^{(l)})dN^{(s)}(u)\Big]\nonumber\\
    &\stackrel{(c)}{=}\sum_{i=1}^s\int_0^\tau E_{F(u)}\Big[w_i\dot{f}_{\boldsymbol{\beta}}(\boldsymbol{x}^{(i)})dN^{(s)}(u)\Big]\nonumber\\ 
    &- \sum_{l=1}^{s}\int_0^\tau E_{F(u)}\Big[w_l\dot{f}_{\boldsymbol{\beta}}(\boldsymbol{x}^{(l)})dN^{(s)}(u)\Big]\nonumber\\
    &= 0
    \label{A3}
\end{align}
where (a) follows from the conditional expectation given the filtration $F(u)$; (b) follows from the fact that $w_l$ and $\dot{f}_{\boldsymbol{\beta^*}}(\boldsymbol{x}^{(l)})$ are known given $F(u)$; (c) follows from the fact that  $E_{\boldsymbol{\beta^*}|F(u)}[dN_i(u))]=w_idN(u)$. Since $E_{\boldsymbol{\beta^*}}[U^{s}(\boldsymbol{\beta^*})]=0$, the parameter $\boldsymbol{\beta}^{(s)}$ is Fisher consistent. In the following, we present a sufficient condition under which such a parameter is a global minmizer of $E[-pl^{(s)}(\boldsymbol{\beta}|\mathcal{D}^{(s)})]$. \\

\textbf{\textit{Corollary 1:}} $\boldsymbol{\beta^{(s)}}$ is a global minimizer of $E[-pl^{(s)}(\boldsymbol{\beta}|\mathcal{D}^{(s)})]$ if $f_{\boldsymbol{\beta}}(\boldsymbol{x}^{(i)})$ is an affine function of $\boldsymbol{\beta}$.\\
\textit{Proof}. Suppose that we choose $f_{\boldsymbol{\beta}}(\boldsymbol{x}^{(l)})$ a convex function of $\boldsymbol{\beta}$. Then, the first term inside summation of $-pl^{(s)}(\boldsymbol{\beta}|\mathcal{D}^{(s)})$ (i.e., $-f_{\boldsymbol{\beta}}(\boldsymbol{x}^{(l)})$) is a concave function. We show that the second term inside summation of $-pl^{(s)}(\boldsymbol{\beta}|\mathcal{D}^{(s)})$, i.e., $log\Big(\sum_{l=1}^s Y_l(u)exp(f_{\boldsymbol{\beta}}(\boldsymbol{x}^{(l)}))\Big)$ is a convex function through the following steps:\\

\textit{Step 1:} $exp(f_{\boldsymbol{\beta}}(\boldsymbol{x}^{(l)}))$ is a convex function because $exp(.)$ is a non-decreasing convex function and $f_{\boldsymbol{\beta}}(\boldsymbol{x}^{(l)})$ is convex \citep[sec.~3.2.4]{Boyd2004}.\\

\textit{Step 2:} $\sum_{l=1}^s Y_l(u)exp(f_{\boldsymbol{\beta}}(\boldsymbol{x}^{(l)}))$ is a convex function because a non-negative weighted sum of convex functions is convex \citep[sec.~3.2.1]{Boyd2004}.\\

\textit{Step 3:} $log\Big(\sum_{l=1}^s Y_l(u)exp(f_{\boldsymbol{\beta}}(\boldsymbol{x}^{(l)}))\Big)$ is a convex function because $log(.)$ is non-decreasing concave function and $\sum_{l=1}^s Y_l(u)exp(f_{\boldsymbol{\beta}}(\boldsymbol{x}^{(l)})$ is a convex function \citep[sec.~3.2.4]{Boyd2004}.\\

\noindent
Therefore, expression $-f_{\boldsymbol{\beta}}(\boldsymbol{x}^{(i)})Y_i(u)+log\Big(\sum_{l=1}^s Y_l(u)exp(f_{\boldsymbol{\beta}}(\boldsymbol{x}^{(l)}))\Big)$ is the sum of a concave function and a convex function. A sufficient condition for convexity of this expression is convexity of $-f_{\boldsymbol{\beta}}(\boldsymbol{x}^{(i)})$ or equivalently concavity of $f_{\boldsymbol{\beta}}(\boldsymbol{x}^{(i)})$. This means that $f_{\boldsymbol{\beta}}(\boldsymbol{x}^{(i)})$ needs to be both convex and concave at the same time. The only functions that satisfy this condition are affine functions \citep[sec.~3.1.1]{Boyd2004}. By choosing $f_{\boldsymbol{\beta}}(\boldsymbol{x}^{(i)})$ as an affine function, $-pl^{(s)}(\boldsymbol{\beta}|\mathcal{D}^{(s)})$ and hence $E[-pl^{(s)}(\boldsymbol{\beta}|\mathcal{D}^{(s)})]$ become convex functions. Therefore, the parameter $\boldsymbol{\beta^{(s)}}$ becomes a global minimizer of $E[-pl^{(s)}(\boldsymbol{\beta}|\mathcal{D}^{(s)})]$. $\square$

Having the loss function $E[-pl^{(s)}(\boldsymbol{\beta}|\mathcal{D}^{(s)})]$ a convex function motivates us to explore the convergence rate of SGD-based minimization algorithms in the next section.

\section{Convergence rate of SGD-based estimate}
\label{secAC}
In this appendix, for the sake of completeness, we prove that our SGD-based estimate can achieve the convergence rate of $O(n^{-1})$ for $f_{\boldsymbol{\beta}}(\bold{x})=\boldsymbol{\beta}^T\bold{x}$, $s=2$, and no ties. We also assume that our covariates are bounded, i.e., there exists some $C < \infty$ such that $||\bold{x}||<C$ with probability $1$ and that we consider a domain of optimization $\bold{B}$ such that $\underset{\boldsymbol{\beta}\in \bold{B}}{\operatorname{max}}\Big\{||\boldsymbol{\beta}^T\boldsymbol{x}||\Big\}<\infty$. We believe these results would hold for $s>2$, however the calculation becomes messier\\   

For the sake of simplicity of proof, we rewrite our objective function as 
\begin{align}
\boldsymbol{\beta}^{(s)} &= \underset{\boldsymbol{\beta}}{\operatorname{argmin}}\Big\{\mathbb{E}_{s}[-pl^{(s)}(\boldsymbol{\beta}|\mathcal{D}^{(s)})]\Big\}\nonumber\\
&= \underset{\boldsymbol{\beta}}{\operatorname{argmin}}\Big\{\mathbb{E}_{s}[L^{(s)}(\boldsymbol{\beta})]\Big\}\nonumber\\
&= \underset{\boldsymbol{\beta}}{\operatorname{argmin}}\Big\{L(\boldsymbol{\beta})\Big\},
\label{B1}
\end{align}
and consider estimating $\boldsymbol{\beta}^{(s)}$ iteratively using SGD from strata of size $s$ through
\begin{align}
\hat{\boldsymbol{\beta}}(m) = \hat{\boldsymbol{\beta}}(m-1) - \gamma_m \times\nabla_{\boldsymbol{\beta}}L^{(s)}(\hat{\boldsymbol{\beta}}(m-1)|\mathcal{D}^{(s)}_m).
\label{B2}
\end{align} 
Authors in \citep{Eric2011} showed that if the loss function $L^{(s)}(\boldsymbol{\beta})$ satisfies the 4 following conditions, then, we can achieve the optimal convergence rate of $O(m^{-1})$ by choosing $\gamma_m=Cm^{-1}$ for the single SGD-based iterates $\hat{\boldsymbol{\beta}}(m)$ and $\gamma_m=C/\sqrt{m}$, $r \in [0.5, 1)$ for averaging over iterates $\tilde{\boldsymbol{\beta}}(m)$ (i.e., Polyak-Ruppert average).\\

\textbf{Condition 1:} Gradient of $L^{(s)}(\boldsymbol{\beta})$ is an unbiased estimate of the gradient $L(\boldsymbol{\beta})$,\\

\textbf{Condition 2:} Gradient of $L^{(s)}(\boldsymbol{\beta})$ is D-Lipschitz-continuous, i.e., $\forall$ $\boldsymbol{\beta}_1, \boldsymbol{\beta}_2 \in \mathbb{R}^p$, there exists $D\geq0$ such that,
\begin{align}
    ||\nabla_{\boldsymbol{\beta}}L^{(s)}(\boldsymbol{\beta}_1)-\nabla_{\boldsymbol{\beta}}L^{(s)}(\boldsymbol{\beta}_2)||\leq D||\boldsymbol{\beta}_1-\boldsymbol{\beta}_2||,
    \label{B3}
\end{align}

\textbf{Condition 3:} $L(\boldsymbol{\beta})=\mathbb{E}_{s}[L^{(s)}(\boldsymbol{\beta})]$ is $\mu$-strongly convex, i.e., $\forall$ $\boldsymbol{\beta}_1, \boldsymbol{\beta}_2 \in \mathbb{R}^p$, there exists $\mu>0$ such that,
\begin{align}
    L(\boldsymbol{\beta}_1) \geq L(\boldsymbol{\beta}_2)+\nabla_{\boldsymbol{\beta}}L(\boldsymbol{\beta}_2)^T(\boldsymbol{\beta}_1-\boldsymbol{\beta}_2)+\frac{\mu}{2}||\boldsymbol{\beta}_1-\boldsymbol{\beta}_2||^2,
    \label{B4}
\end{align}

\textbf{Condition 4:} Variance of the gradient of $L^{(s)}(\boldsymbol{\beta})$ is bounded, i.e., there exists $\sigma^2 \in \mathbb{R}_{+}$ such that
\begin{align}
    \mathbb{E}(||\nabla_{\boldsymbol{\beta}}L^{(s)}(\boldsymbol{\beta}^*)||^2) \leq \sigma^2, w.p.1,
    \label{B5}
\end{align}

\noindent
In the following, we show that the loss function in our framework, $L^{(s)}(\boldsymbol{\beta})$ satisfies all four conditions above.\\

\textbf{Proof of Condition 1:} This condition is automatically satisfied based on the definition of $L(\boldsymbol{\beta})=\mathbb{E}_{s}[L^{(s)}(\boldsymbol{\beta})]$ in \eqref{B1} and that $\nabla_{\boldsymbol{\beta}}L(\boldsymbol{\beta})=\nabla_{\boldsymbol{\beta}}\mathbb{E}_{s}[L^{(s)}(\boldsymbol{\beta})]=\mathbb{E}_{s}[\nabla_{\boldsymbol{\beta}}L^{(s)}(\boldsymbol{\beta})]$.\\

\textbf{Proof of condition 2:} The loss function $L^{(s)}(\boldsymbol{\beta})$ belongs to $C^\infty$ continuous function family and proving \eqref{B3} is equivalent to proving
\begin{align}
    \exists\, D\geq0,\, s.t., \,\forall \nu \in S_\nu=\{\nu: ||\nu||_2 = 1\},\quad \nu^T\nabla^2_{\boldsymbol{\beta}}L^{(s)}(\boldsymbol{\beta})\nu \leq D
    \label{B6}.
\end{align}
For $f_{\boldsymbol{\beta}}(\boldsymbol{x})=\boldsymbol{\beta}^T\boldsymbol{x}$, $s=2$ and assuming no ties, $\nabla^2_{\boldsymbol{\beta}}L^{(s)}(\boldsymbol{\beta})$ can be simplified as 
\begin{align}
    \nabla^2_{\beta}\Big\{L^{(s)}(\boldsymbol{\beta})\Big\}&= w(1-w)\bold{X}\begin{bmatrix}
1 & -1 \\
-1 & 1
\end{bmatrix}\bold{X}^T\nonumber\\
&\times\Big(1(\delta_1=1)1(y_1<y_2)+1(\delta_2=1)1(y_2<y_1)\Big)\nonumber\\
&=w(1-w)(\bold{x}^{(1)}-\bold{x}^{(2)})(\bold{x}^{(1)}-\bold{x}^{(2)})^T\nonumber\\
&\times\Big(1(\delta_1=1)1(y_1<y_2)+1(\delta_2=1)1(y_2<y_1)\Big)
    \label{B7}
\end{align}
where $\bold{X} = [\bold{x}^{(1)}, \bold{x}^{(2)}]$; $w = w_1(\bold{x}^{(1)},\bold{x}^{(2)}, \boldsymbol{\beta})= 1-w_2(\bold{x}^{(1)},\bold{x}^{(2)}, \boldsymbol{\beta})=\frac{\exp{(\boldsymbol{\beta}^T\boldsymbol{x}^{(1)})}}{\exp{(\boldsymbol{\beta}^T\boldsymbol{x}^{(1)})}+\exp{(\boldsymbol{\beta}^T\boldsymbol{x}^{(2)})}}$. Note that $1(\delta_1=1)1(y_1<y_2)+1(\delta_2=1)1(y_2<y_1) \leq 1$ and that $w(1-w)\leq0.25$ because $0\leq w \leq 1$. Therefore we have
\begin{align}
    \nu^T\nabla^2_{\beta}\Big\{L^{(s)}(\boldsymbol{\beta})\Big\}\nu&\leq 0.25 \times \nu^T(\bold{x}^{(1)}-\bold{x}^{(2)})(\bold{x}^{(1)}-\bold{x}^{(2)})^T\nu\nonumber\\
    &=0.25\times ||\nu^T(\bold{x}^{(1)}-\bold{x}^{(2)})||_2^2\nonumber\\
    & \leq 0.25\times ||(\bold{x}^{(1)}-\bold{x}^{(2)})||_2^2\times ||\nu||_2^2\nonumber\\
    & \stackrel{(a)}{\leq} 0.25 \times (||\bold{x}^{(1)}||_2^2 + ||\bold{x}^{(2)}||_2^2)\nonumber\\
    & \leq 0.25 \times 2 \times max(||\bold{x}^{(1)}||_2^2, ||\bold{x}^{(2)}||_2^2)\nonumber\\
    & \leq 0.5 \times {\underset{i}{\operatorname{max}}{||\bold{x}^{(i)}||_2^2}},
    \label{B8}
\end{align}
where $(a)$ follows the triangle inequality. Therefore, by assuming our covariates $\boldsymbol{x}$ are bounded, the gradient of $L^{(s)}(\boldsymbol{\beta})$ is $D-$Lipschitz-continuous with $D=0.5 \times {\underset{i}{\operatorname{max}}{||\bold{x}^{(i)}||^2}}$. This completes the proof of \textbf{Condition 2}. $\square$\\

\textbf{Proof of condition 3:} The loss function $L(\boldsymbol{\beta})$ belongs to the $C^\infty$ continuous function family and proving \eqref{B4} is equivalent to proving
\begin{align}
    \exists\, \mu>0, \, s.t.\, \forall \nu \in S_\nu=\{\nu: ||\nu||_2 = 1\},\quad \nu^TI(\boldsymbol{\beta})\nu \geq \mu
    \label{B9}.
\end{align}
where $I(\boldsymbol{\beta}) = \mathbb{E}_{s}[\nabla^2_{\boldsymbol{\beta}}L^{(s)}(\boldsymbol{\beta})]$ is the expected Hessian matrix. Starting from \eqref{B7}, $I(\boldsymbol{\beta})$ can be written as
\begin{align}
    I(\boldsymbol{\beta}) &=E_s\Big[\nabla^2_{\beta}\Big\{L^{(s)}(\boldsymbol{\beta})\Big\}\Big]\nonumber\\
    &=E_{X, Y, \Delta}\Big[\nabla^2_{\beta}\Big\{L^{(s)}(\boldsymbol{\beta})\Big\}\Big]\nonumber\\
    &\stackrel{(a)}{=}E_{X}\Bigg[E_{Y|X}\Big[E_{\Delta|X, Y}[\nabla^2_{\beta}\Big\{L^{(s)}(\boldsymbol{\beta})\Big\}]\Big]\Bigg]\nonumber\\
    &\stackrel{(b)}{=}(1-p_c)E_{X}\Big[w(1-w)(\bold{x}^{(1)}-\bold{x}^{(2)})^T(\bold{x}^{(1)}-\bold{x}^{(2)})\Big]\nonumber\\
    &=(1-p_c) \int_{X^{(1)}, X^{(2)}}\Big[w(1-w)(\bold{x}^{(1)}-\bold{x}^{(2)})^T(\bold{x}^{(1)}-\bold{x}^{(2)})\Big]P_{X}(x)\nonumber\\
    &\stackrel{(c)}{=} (1-p_c) \int_{\bold{Z}}w(1-w)\bold{Z}^T\bold{Z}P_{\bold{Z}}(\bold{z})
    \label{B10}
\end{align}
where $(a)$ follows the expansion of the intersection using conditional probabilities; $(b)$ follow from the fact that $p_c=E_{\Delta}[1(\delta_i=0))]=P_{\Delta}(\delta_i=0)$ is the probability of censoring and that we have $E_{Y}[1(y_1<y_2))+1(y_2<y_1)] = Pr(y_1<y_2) + Pr(y_2<y_1)=1$; $(c)$ follows from the change of variable $\bold{Z}=\bold{X}^{(1)}-\bold{X}^{(2)}$. Note that since random variables $\bold{X}^{(1)}$ and $\bold{X}^{(2)}$ have density with respect to the Lebesgue measure, random variable $\bold{Z}$ also has a density with respect to the Lebesgue measure. Then we can write $\nu^TI(\boldsymbol{\beta})\nu$ as
\begin{align}
    \nu^TI(\boldsymbol{\beta})\nu &= (1-p_c)\int_{\bold{Z}}w(1-w)\nu^T\bold{Z}^T\bold{Z}\nu P_{\bold{Z}}(\bold{z})\nonumber\\
    &= (1-p_c)\int_{\bold{Z}}w(1-w)||{\nu}^T\bold{z}||_2^2P_{\bold{Z}}(\bold{z})
    \label{B11}
\end{align}

Now we prove strong convexity of $L(\boldsymbol{\beta})$ by contradiction. The negation of statement \eqref{B9} is:
\begin{align}
    \forall \mu >0,\, \exists \nu_\mu \in S_\nu=\{\nu: ||\nu||_2 = 1\},\, s.t.\quad \nu_\mu^TI(\boldsymbol{\beta})\nu_\mu < \mu.
    \label{B12}
\end{align}

\textit{Claim 1:} Suppose the statement in \ref{B12} holds (or equivalently \eqref{B9} does not hold), then there exists a $\nu^* \in S_\nu=\{\nu: ||\nu||_2 = 1\}$ such that we have $P_\bold{z}(||{\nu^*}^T\bold{z}||_2=0)=1$.  \\

\textit{Proof:} Since we assumed $\boldsymbol{\beta}^T\boldsymbol{x}$ is bounded, there exists a constant $0<C_w<1$ such that $C_w<w<1-C_w$, and hence $w(1-w)>C_w^2$. Therefore, using \eqref{B11}, the statement \eqref{B12} implies that for any $\mu>0$ there exists $\nu_\mu \in S_\nu$ such that
\begin{align}
    \mu &> \nu_\mu^TI(\boldsymbol{\beta})\nu_\mu=(1-p_c)\int_{\bold{Z}}w(1-w)||\nu_\mu^T\bold{z}||_2^2P_{\bold{Z}}(\bold{z})\nonumber\\
    &\stackrel{(a)}{>}(1-p_c)C_w^2\int_{\bold{Z}}||\nu_\mu^T\bold{z}||_2^2P_{\bold{Z}}(\bold{z})\nonumber\\
    &=(1-p_c)C_w^2\int_{\bold{Z}}\Big(1(||\nu_\mu^T\bold{z}||_2^2<\epsilon))+1(||\nu_\mu^T\bold{z}||_2^2>\epsilon)\Big)||\nu_\mu^T\bold{z}||_2^2P_{\bold{Z}}(\bold{z})\nonumber\\
    &\geq (1-p_c)C_w^2\int_{\bold{Z}}1(||\nu_\mu^T\bold{z}||_2^2>\epsilon)||\nu_\mu^T\bold{z}||_2^2P_{\bold{Z}}(\bold{z})\nonumber\\
    &> (1-p_c)C_w^2\epsilon P_{\bold{Z}}(||\nu_\mu^T\bold{z}||_2^2>\epsilon).
    \label{B13}
\end{align}
Note that $\epsilon$ is an arbitrary positive value. Therefore, for such a $\nu_\mu$, expression \eqref{B13} is equivalent to
\begin{align}
    P_{\bold{Z}}(||\nu_\mu^T\bold{z}||_2^2>\epsilon) < \frac{\mu}{(1-p_c)C_w^2\epsilon}\quad \forall \mu, \, \epsilon>0
    \label{B14}
\end{align}
or equivalently,
\begin{align}
    P_{\bold{Z}}(||\nu_\mu^T\bold{z}||_2^2<\epsilon) >1- \frac{\mu}{(1-p_c)C_w^2\epsilon}\quad \forall \mu, \, \epsilon>0.
    \label{B15}
\end{align}
 Thus, there exists an infinite sequence $\nu_1, \dots, \nu_k, \dots$ such that for the choices of $\mu_k(\delta)=\delta(1-p_c)C_w^2\epsilon$ and $\epsilon_k=\frac{1}{k}$ we have
\begin{align}
    \forall\,\delta>0\,, \exists\,K>0\,, s.t.\,, \forall\, k > K: Pr(||\nu_k^T\bold{z}||_2\leq \frac{1}{k}) > 1-\delta.
    \label{B16}
\end{align}
Since $S_\nu$ is a compact space, this sequence has an infinite subsequence converging to a point $\nu^* \in S_\nu$. For such a converging point $\nu^*$, we can write, for all $k$,
\begin{align}
    ||{\nu^*}^T\bold{z}||_2 &= ||(\nu^*-\nu_{k}+\nu_{k})^T\bold{z}||_2\nonumber\\
    & \stackrel{(a)}{\leq} ||(\nu^*-\nu_{k})^T\bold{z}||_2+||\nu_{k}^T\bold{z}||_2\nonumber\\
    & \leq ||\nu^*-\nu_{k}||_2||\bold{z}||_2+||\nu_{k}^T\bold{z}||_2\nonumber\\
    & \stackrel{(b)}{\leq} ||\nu^*-\nu_{k}||_2C_z+||\nu_{k}^T\bold{z}||_2,
    \label{B17}
\end{align}
where $(a)$ follows the triangle inequality and $(b)$ follows the boundedness of variable $\bold{Z}$ (i.e., $||\bold{z}||_2\leq C_z$) due to boundedness of variable $\bold{X}$. From \eqref{B16}, we may write
\begin{align}
    \forall\,\delta>0\,, \epsilon>0\,, \exists\,K_1>0\,, s.t.\,, \forall\, k > K_1: Pr(||\nu_k^T\bold{z}||_2\leq \frac{\epsilon}{2}) > 1-\delta,
    \label{B18}
\end{align}
and since $\nu_k$ converges to $\nu^*$, we may write
\begin{align}
    \forall\,\delta>0\,, \epsilon>0\,, \exists\,K_2>0\,, s.t.\,, \forall\, k > K_2: Pr(||\nu^*-\nu_k||_2\leq \frac{\epsilon}{2C_z}) > 1-\delta.
    \label{B19}
\end{align}
Then $\forall\,\delta>0$ and $\forall\,\epsilon>0$, for $K_{max}=max(K_1, K_2)$ we have that for all $k>K_{max}$
\begin{align}
    Pr(||{\nu^*}^T\bold{z}||_2 &\leq ||\nu^*-\nu_{k}||_2C_z+||\nu_{k}^T\bold{z}||_2\nonumber\\
    & \leq \frac{\epsilon}{2C_z}\times C_z + \frac{\epsilon}{2}=\epsilon) > 1-\delta,
    \label{B20}
\end{align}
taking $\epsilon, \delta \rightarrow 0$, we have that $P_{\bold{Z}}(||{\nu^*}^T\bold{z}||_2=0)=1$. This completes the proof of \textit{Claim 1}. $\square$ \\

\textit{Claim 2:} If random variable $\bold{X}^{(1)}$ and $\bold{X}^{(2)}$ have density with respect to the Lebesgue measure, then we have $P_{\bold{Z}}(||{\nu}^T\bold{z}||_2=0)<1$ for any $\nu \in S_\nu=\{\nu: ||\nu||_2 = 1\}$.\\

\textit{Proof:} Random variables $\bold{X}^{(1)}$ and $\bold{X}^{(2)}$ have density with respect to the Lebesgue measure. Therefore, variable $\bold{Z}$ also has a density with respect to the Lebesgue measure. Therefore, for any $\nu \in S_\nu=\{\nu: ||\nu||_2 = 1\}$, $\bold{Z}$ cannot only be on a plane orthogonal to $\nu$. In other words, 
\begin{align}
    \forall \nu \in S_\nu=\{\nu: ||\nu||_2 = 1\}, \exists \epsilon_{\nu}>0\, \text{and}\, \exists \delta_{\nu}>0\, s.t.\quad P_{\bold{Z}}(||{\nu^*}^T\bold{z}||_2^2>\epsilon_{\nu}) > \delta_{\nu}.
    \label{B21}
\end{align}
Then for any $\nu \in S_\nu=\{\nu: ||\nu||_2 = 1\}$, the integral $\int_{\bold{Z}}||\nu^T\bold{z}||_2^2dP_{\bold{Z}}(\bold{z})$ can be written as
\begin{align}
    \int_{\bold{Z}}||\nu^T\bold{z}||_2^2dP_{\bold{Z}}(\bold{z})&=\int_{\bold{Z}}\Big(1(||\nu^T\bold{z}||_2^2<\epsilon_{\nu}))+1(||\nu^T\bold{z}||_2^2>\epsilon_{\nu})\Big)||\nu^T\bold{z}||_2^2dP_{\bold{Z}}(\bold{z})\nonumber\\
    &\geq \int_{\bold{Z}}1(||\nu^T\bold{z}||_2^2>\epsilon_{\nu})||\nu^T\bold{z}||_2^2dP_{\bold{Z}}(\bold{z})\nonumber\\ 
    &> \epsilon_{\nu}\int_{\bold{Z}}1(||\nu^T\bold{z}||_2^2>\epsilon_{\nu})dP_{\bold{Z}}(\bold{z})\nonumber\\
    &= \epsilon_{\nu} P_{\bold{Z}}(||{\nu^*}^T\bold{z}||_2^2>\epsilon_{\nu})\nonumber\\
    &> \epsilon_{\nu} \delta_{\nu} > 0,
    \label{B22}
\end{align}
indicating that $P_{\bold{Z}}(||{\nu}^T\bold{z}||_2=0)<1$ (or $P_{\bold{Z}}(||{\nu}^T\bold{z}||_2=0)\not=1$). This completes the proof of \textit{Claim 2}. $\square$\\

\textit{Claim 1} indicated that if the strong convexity statement \eqref{B9} does not hold, we must have $P_\bold{z}(||{\nu^*}^T\bold{z}||_2=0)=1$. This is in contradiction with the result of \textit{Claim 2} indicating that for random variable $\bold{Z}$ with Lebesgue density we have $P_\bold{z}(||{\nu^*}^T\bold{z}||_2=0)\not=1$. Therefore, \eqref{B12} is not true and the statement \eqref{B9} holds, i.e., $L(\boldsymbol{\beta})$ is strongly convex. This completes the proof of \textbf{Condition 3}. $\square$ \\

\textbf{Proof of condition 4:} The gradient of $L^{(s)}(\boldsymbol{\beta})$ may be written as
\begin{align}
    \nabla_{\boldsymbol{\beta}}L^{(s)}(\boldsymbol{\beta}) = \sum_{i=1}^s1(\delta_i=1)(\bold{x}^{(i)}-\sum_{j \in R_i}w_j\boldsymbol{x}^{(j)}).
    \label{B23}
\end{align}
Then we can write
\begin{align}
    ||\nabla_{\boldsymbol{\beta}}L^{(s)}(\boldsymbol{\beta})||^2 &= ||\sum_{i=1}^s1(\delta_i=1)(\bold{x}^{(i)}-\sum_{j \in R_i}w_j\boldsymbol{x}^{(j)})||^2 \nonumber\\
    &\stackrel{(a)}{\leq} \Big(\sum_{i=1}^s(||1(\delta_i=1)||\times||\bold{x}^{(i)}-\sum_{j \in R_i}w_j\boldsymbol{x}^{(j)})||\Big)^2\nonumber\\
    &\stackrel{(b)}{\leq} \Big(\sum_{i=1}^s(||\bold{x}^{(i)}||+||\sum_{j \in R_i}w_j\boldsymbol{x}^{(j)}||)\Big)^2\nonumber\\
    &\stackrel{(c)}{\leq} \Big(\sum_{i=1}^s(||\bold{x}^{(i)}||+ {\underset{j \in R_i}{\operatorname{max}}{||\boldsymbol{x}^{(j)}||)}} \Big)^2\nonumber\\
    &\leq 4S^2{\underset{i}{\operatorname{max}}{||\boldsymbol{x}^{(i)}||^2}}
    \label{B24}
\end{align}
where $(a)$ follows from the triangle inequality; $(b)$ and $(c)$ follow from the triangle inequality and that $w_j\leq1$. Therefore, given boundedness of covariates $\bold{x}$, we can choose $\sigma=2S{\underset{i}{\operatorname{max}}{||\boldsymbol{x}^{(i)}}}||$. This completes the proof of \textbf{Condition 4}. $\square$ \\

We showed that all four conditions are satisfied and thus the results in \citep{Eric2011} give us our claimed convergence rates for both single SGD-based estimates and averaging over estimates (Polyak-Ruppert average).

\section{Mini-batches, Moment-based-learning-rate, and Multiple Epochs}
\label{secAD}
For small strata size $s$, the stochastic gradients given in \eqref{e7} are potentially quite noisy. In these cases, rather than using only a single stratum to calculate the stochastic gradient, it may be preferable to average gradients over multiple strata. This is known as a using a mini-batch \citep{Ruder2016} in the SGD literature. In this survival context, we use batches of strata, so a larger strata size $s$ can eliminate the need for mini-batch sizes of greater than $1$.

Choosing a reasonable value for the hyper-parameter $\gamma_m$ (the learning rate) is critical for good practical performance of the algorithm (see Appendix \ref{secAG11} for related simulation results). A number of publications have developed methods for adaptively selecting the learning rate using moments, including {\tt Adagrad} \citep{Duchi2011}, and {\tt ADAM} \citep{Kingma2014}. However, it was shown that {\tt ADAM} may not converge in some settings and an updated version called {\tt AMSGrad} has been proposed \citep{Sashank2019}. We found substantially improved performance on simulated data with {\tt AMSGrad} over the simple non-adaptive updating rule given in \eqref{e8}, especially in combination with averaging over iterates as discussed in Section 3.3.2.

In practice, taking only a single pass over the data leads to poor empirical performance. We generally use multiple passes (or ``epochs''). In particular, for each epoch, we (1) randomly partition our data into strata; then (2) use the last updated iterate of $\hat{\boldsymbol{\beta}}$ from the previous epoch as the initial iterate of $\hat{\boldsymbol{\beta}}$ for the new epoch; and finally (3) apply a full pass of stochastic gradient descent (with e.g., averaging and momentum) over the partitioned data. In practice we found that $\sim 100$ epochs was more than enough for very robust convergence (see the top right panel of Figure~1 in Section 6.2).

\section{Implementation of streaming and non-streaming algorithms using BigSurvSGD}
\label{secAE}
In this section, we present the implementations of both streaming and non-streaming mini-batch stochastic gradient descent algorithms using our proposed framework BigSurvSGD. Without loss of generality, we only present algorithms without the moment-based step-size adaptation. 
\subsection{Implementation of streaming algorithm}
\label{secAE1}
As we discussed in the main manuscript, our proposed framework facilitates the implementation of an algorithm using SGD that handles the streaming data. This implementation is very straightforward: We update the estimate (i.e., $\hat{\boldsymbol{\beta}}$) in a streaming fashion using each new coming mini-batch of data. Therefore, our algorithm does not run into the memory issues because we do not need to collect all the data before the estimation (as would be required by {\tt coxph()}). Algorithm \ref{algStreaming} summarizes the implementation of a streaming algorithm using our proposed framework where we use a single stratum for each step (i.e., mini-batches of size 1). Extension to mini-batches of size greater than 1 is very straightforward (see next section).\\

\begin{algorithm}[H]
\SetAlgoLined
\KwResult{$\hat{\boldsymbol{\beta}}$}
\textbf{Initialization}:\\
\quad Choose strata size $s$\\
\quad $\hat{\boldsymbol{\beta}}(0)=\boldsymbol{0}$\\
\quad Choose $C$ for $\gamma_m=\frac{C}{\sqrt{m}}$

\For{($m=1, 2, ...$)}{
Draw data $\mathcal{D}^{(s)}_m$ including $s$ patients from patient $s\times (m-1)+1$ to patient $s\times m$\\
Update the estimate $\hat{\boldsymbol{\beta}}(m)$ using
\begin{align*}
\hat{\boldsymbol{\beta}}(m) = \hat{\boldsymbol{\beta}}(m-1) + \gamma_m \times \dot{pl}^{(s)}(\hat{\boldsymbol{\beta}}(m-1)|\mathcal{D}^{(s)}_m)
\end{align*}
}
\textbf{Output}:\\
Calculate averaged over iterate estimates as $\widetilde{\boldsymbol{\beta}}(m)=\frac{1}{m}\sum_{l=1}^{m}\widehat{\boldsymbol{\beta}}(l)$, $l=1, 2, \dots, m$
\caption{Implementation of streaming algorithm using BigSurvSGD}
\label{algStreaming}
\end{algorithm}

\subsection{Implementation of a non-streaming mini-batch stochastic gradient descent algorithm}
\label{secAE2}
In many cases, it is of interest to use each datum more than once (this can empirically improve performance). In these cases, we still use a mini-batch stochastic gradient descent algorithm. We split our data evenly into mini-batches, each with $K$ strata of size $s$. Then we iteratively update the estimate using the batches of strata. Using multiple, rather than single strata per batch results in more stable (less noisy) gradients. We use learning rate $\frac{\gamma_m}{K}$ instead of $\gamma_m$ since the gradient in each step of mini-batch gradient descent is the sum of gradients over $K$ strata with size $s$. One strength of this implementation compared to {\tt coxph()} is we can read the batches of strata chunk-by-chunk from the hard drive. Therefore, this implementation handles large amounts of data without running into the memory issues. Algorithm \ref{algMiniBatch} summarizes the implementation of a non-streaming mini-batch stochastic gradient descent algorithm using our proposed framework. \\

\begin{algorithm}[H]
\SetAlgoLined
\KwResult{$\hat{\boldsymbol{\beta}}$}
\textbf{Initialization}:\\
\quad Choose strata size $s$\\
\quad Choose batch size $K$\\
\quad Choose number of epochs $n_E$\\
\quad $m=0$\\
\quad $\hat{\boldsymbol{\beta}}(0)=\boldsymbol{0}$\\
\quad Choose $C$ for $\gamma_m=\frac{C}{\sqrt{m}}$\\
\For{($n_e=1, 2, ..., n_E$)}{
Divide data randomly into $n_B$ disjoint batches, $b=1, 2, \dots, n_B$. Each batch includes $K$ disjoint strata of patients, $\mathcal{D}^{(s)}_{b,1}$, ..., $\mathcal{D}^{(s)}_{b,K}$ with size $s$.\\

\For{($b=1, 2, ..., n_B$)}{
$m = m+1$\\
Update the estimate of $\hat{\boldsymbol{\beta}}(m)$ using
\begin{align*}
\hat{\boldsymbol{\beta}}(m) = \hat{\boldsymbol{\beta}}(m-1) + \frac{\gamma_m}{K} \times \sum_{k=1}^{K}\dot{pl}^{(s)}(\hat{\boldsymbol{\beta}}(m-1)|\mathcal{D}_{b,k}^{(s)})
\end{align*}
}
}
\textbf{Output}:\\
Calculate averaged over iterate estimates as $\widetilde{\boldsymbol{\beta}}(m)=\frac{1}{m}\sum_{l=1}^{m}\widehat{\boldsymbol{\beta}}(l)$, $l=1, 2, \dots, m$
\caption{Non-streaming mini-batch gradient descent algorithm using BigSurvSGD}
\label{algMiniBatch}
\end{algorithm}

\section{Implementation of algorithms for constructing confidence interval}
\label{secAF}
In this section, we present the implementations of both the plug-in and bootstrap algorithms for constructing a confidence interval.
\subsection{Constructing a $(1-\alpha)100\%$ confidence interval using plug-in method}
\label{secAF1}
Algorithm~\ref{algCIplug} presents details for constructing a $(1-\alpha)100\%$ confidence interval using the plug-in approach proposed in Section~5.1 of the main manuscript. To reduce the computational complexity, we use $n_o$ ($n_o << {n-1\choose s-1}$) strata sample per observation to estimate the standard error. In practice, sampling $100\sim1000$ strata per observation suffices to get a near nominal coverage (see Appendix~\ref{secAE2} for related results). After estimating the standard error, we construct the $(1-\alpha)100\%$ confidence interval following the asymptotic normality of estimator $\widetilde{\boldsymbol{\beta}}$.  \\

\begin{algorithm}[!tp]
\SetAlgoLined
\KwResult{$\widetilde{\boldsymbol{\beta}}$ and $(1-\alpha)100\%$}
\textbf{Initialization}:\\
\quad Choose strata size $s$\\
\quad Choose $C$ for $\gamma_m=\frac{C}{\sqrt{m}}$\\
\quad Specify significance level $\alpha$\\
\quad Calculate $\widetilde{\boldsymbol{\beta}}$ using Algorithm~\ref{algMiniBatch}\\
\quad Specify $n_o$ ($n_o << {n-1\choose s-1}$), number of strata sample per observation to estimate the standard error\\

\For{($i=1, 2, ...n$)}{
Randomly choose $\mathcal{M}$, a set of $M$ distinct subsets of $s-1$ elements from $\{1, 2, \dots, n\}\backslash \{i\}$. Then estimate $\widehat{r}(\mathcal{D}_i, \widetilde{\boldsymbol{\beta}})$ using
\begin{align*}
    \widehat{r}(\mathcal{D}_i, \widetilde{\boldsymbol{\beta}}) = \frac{1}{M} \sum_{\{i_2, i_3, \dots, i_s\} \in \mathcal{M}}\frac{\partial}{\partial \boldsymbol{\beta}}\{-pl^{(s)}(\boldsymbol{\beta}|\mathcal{D}^{(s)}=\{\mathcal{D}_i, \mathcal{D}_{i_2}, \dots, \mathcal{D}_{i_s}\})\}|_{\boldsymbol{\beta}=\widetilde{\boldsymbol{\beta}}}
\end{align*}
}
Estimate $\widehat{\boldsymbol{V}}_n$ and $\widehat{\boldsymbol{H}}_n$ using
\begin{align*}
    \widehat{\boldsymbol{V}}_n &= \frac{s^2}{n} \sum_{i=1}^{n}\widehat{r}(\mathcal{D}_i, \widetilde{\boldsymbol{\beta}})\widehat{r}(\mathcal{D}_i, \widetilde{\boldsymbol{\beta}})^T \nonumber\\
    \widehat{\boldsymbol{H}}_n &= \left[\frac{\partial}{\partial \boldsymbol{\beta}} \frac{1}{n} \sum_{i=1}^{n} \widehat{r}(\mathcal{D}_i, \boldsymbol{\beta})\right]_{\boldsymbol{\beta}=\widetilde{\boldsymbol{\beta}}}
\end{align*}

\textbf{Output}:\\ 
Estimate standard error of $\widetilde{\boldsymbol{\beta}}$ using $\widehat{\overrightarrow{SE}}=\sqrt{diag(\widehat{\bold{H}}_n^{-1}\bold{\widehat{V}}_n\bold{\widehat{H}}_n^{-1})/n}$. Report the $(1-\alpha)100\%$ CI as $(\widetilde{\boldsymbol{\beta}}-z_{\frac{\alpha}{2}}\times \widehat{\overrightarrow{SE}}, \widetilde{\boldsymbol{\beta}}+z_{\frac{\alpha}{2}} \times \widehat{\overrightarrow{SE}})$ where $z_{\frac{\alpha}{2}}$ is the upper $\frac{\alpha}{2}$ critical value for the standard normal distribution.\\
\caption{Construction of a $(1-\alpha)100\%$ CI using plug-in method}
\label{algCIplug}
\end{algorithm} 

\subsection{Constructing a $(1-\alpha)100\%$ confidence interval using bootstrap method}
\label{secAF2}
Algorithm~\ref{algCIboot} presents details for constructing a $(1-\alpha)100\%$ confidence interval using the bootstrap approach proposed in Section 5.2 of the main manuscript. The first step of this algorithm is to get an initial estimate  $\widetilde{\boldsymbol{\beta}}$ using our proposed Algorithm \ref{algMiniBatch} in Appendix \ref{secAE2}. The bootstrap steps for getting the bootstrapped estimates $\widetilde{\boldsymbol{\beta}}^b$, $b=1, 2, \dots, B$ are exactly the same as what we presented in Algorithm \ref{algMiniBatch} except we use the bootstrap resamples with initialization $\widehat{\boldsymbol{\beta}}^b(0)=\widetilde{\boldsymbol{\beta}}, b=1, 2, \dots, B$. Then, following the fact that $\sqrt{n}(\widetilde{\boldsymbol{\beta}} - \boldsymbol{\beta})$ and $\sqrt{n}(\widetilde{\boldsymbol{\beta}}^b - \widetilde{\boldsymbol{\beta}})$ converges in distribution to the same limiting distribution, we construct the $(1-\alpha)100\%$ for $\widetilde{\boldsymbol{\beta}}$. \\

\begin{algorithm}[H]
\SetAlgoLined
\KwResult{$(1-\alpha)100\%$ confidence interval}
\textbf{Initialization}:\\
\quad Choose strata size $s$\\
\quad Choose $C$ for $\gamma_m=\frac{C}{\sqrt{m}}$\\
\quad Specify significance level $\alpha$\\
\quad Calculate $\widetilde{\boldsymbol{\beta}}$ using Algorithm~\ref{algMiniBatch}\\

\For{b=1, 2, \dots, B}{
Resample $\mathcal{D}^{(n)}$ with replacement to get the $b$-th bootstrap resample $\mathcal{D}^{(n),b}$\\
Use Algorithm \ref{algMiniBatch} to calculate $\widetilde{\boldsymbol{\beta}}^b$ with resampled data $\mathcal{D}^{(n),b}$, and initialization $\widetilde{\boldsymbol{\beta}}^b(0)=\widetilde{\boldsymbol{\beta}}$
}
\textbf{Output}:\\ 
Calculate $q_{\frac{\alpha}{2}}$ and $q_{1-\frac{\alpha}{2}}$ as the $\frac{\alpha}{2}$ and $1-\frac{d}{2}$ quantiles of $\widetilde{\boldsymbol{\beta}}^b - \widetilde{\boldsymbol{\beta}}$, $b=1, 2, \dots, B$. Report the $(1-\alpha)100\%$ CI as $(\widetilde{\boldsymbol{\beta}}-q_{1-\frac{\alpha}{2}}, \widetilde{\boldsymbol{\beta}}-q_{\frac{\alpha}{2}})$.
\caption{Estimate $(1-\alpha)100\%$ confidence interval using bootstrap method}
\label{algCIboot}
\end{algorithm}


\section{Complementary results}
\label{secAG}
\subsection{Simulated data}
\label{secAG1}
\subsubsection{Effect of misspecification of the initial learning rate}
\label{secAG11}
\begin{figure}[!tb]
    \centering
    \begin{subfigure}[b]{0.496\textwidth}
            \centering
            \includegraphics[width=\textwidth]{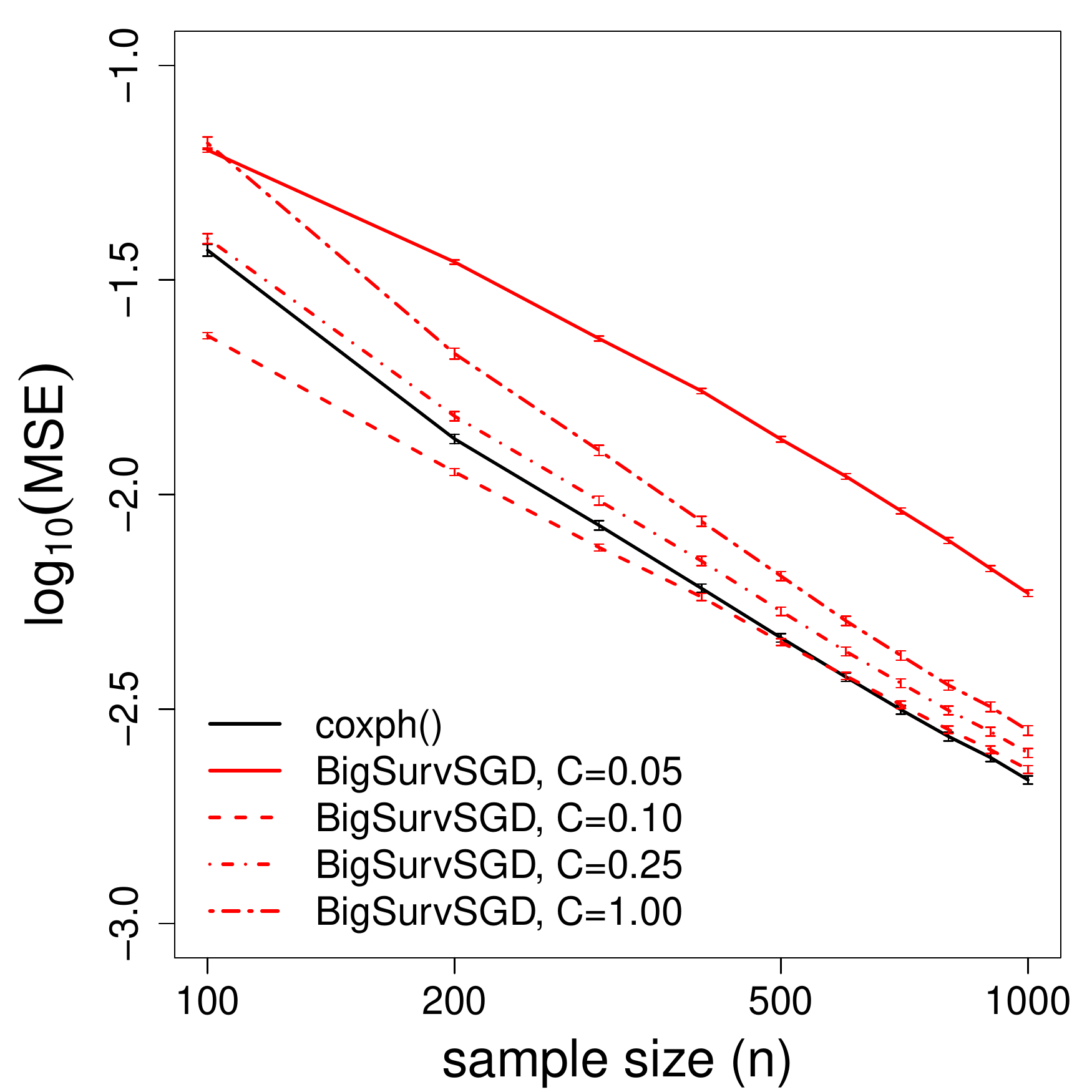}
    \end{subfigure}
    \begin{subfigure}[b]{0.496\textwidth}  
        \centering 
        \includegraphics[width=\textwidth]{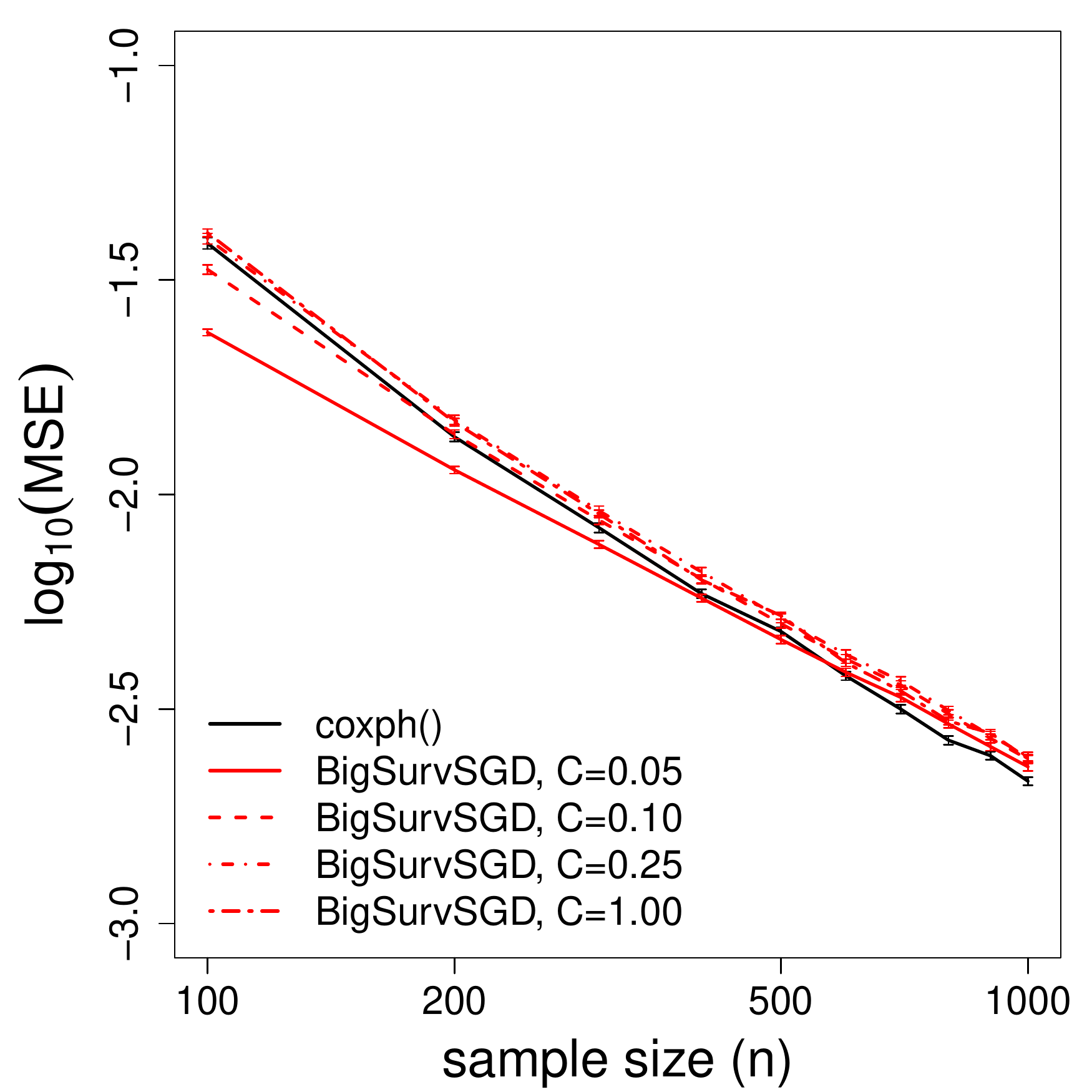}
    \end{subfigure}
    \caption[ The average and standard deviation of critical parameters ]
    {\small (left panel) $\operatorname{\log}_{10}$(MSE) for varying values of $C$ with 100 epochs; (right panel) $\operatorname{\log}_{10}$(MSE) for varying values of $C$ with 1000 epochs. We choose mini-batch size $K=1$, strata size $S=20$, $P=10$ features, and probability of censoring $p_c=0.2$.} 
    \label{fig04}
\end{figure}

Choosing an appropriate learning rate may help improve the statistical efficiency of our estimator for a fixed number of epochs. In practice, we have found that {\tt AMSGrad} is much more robust to misspecification of the learning rate. The left panel in Figure~\ref{fig04} compares MSE with varying choices of $C$ for the learning rate defined as $\gamma_m=\frac{C}{\sqrt{m}}$. We use AveAMSGrad over 1000 simulated datasets (based on data generation mechanism explained in Section 6.1 with probability of censoring $p_c=0.2$) with 100 and 1000 epochs. We observe choosing an appropriate value of $C$ (around 0.1) gives strong performance. Choosing smaller or larger values of $C$ may result in the statistical inefficiency for $100$ epochs: Smaller values of $C$ (e.g., $C=0.05$) result in too little "learning" in each step of the SGD algorithm, while higher values of $C$ result in too much noise. The right panel of Figure~\ref{fig04} considers the same settings as the left panel except we increase the number of epochs from 100 to 1000. As we see, our algorithm with higher number of epochs is relatively robust to a wider range of $C$ around the empirically optimal value.

\subsubsection{Effect of number of sample strata on coverage for plugin approach}
\label{secAG12}
\begin{figure}[!tb]
\centering
\includegraphics[scale=0.7]{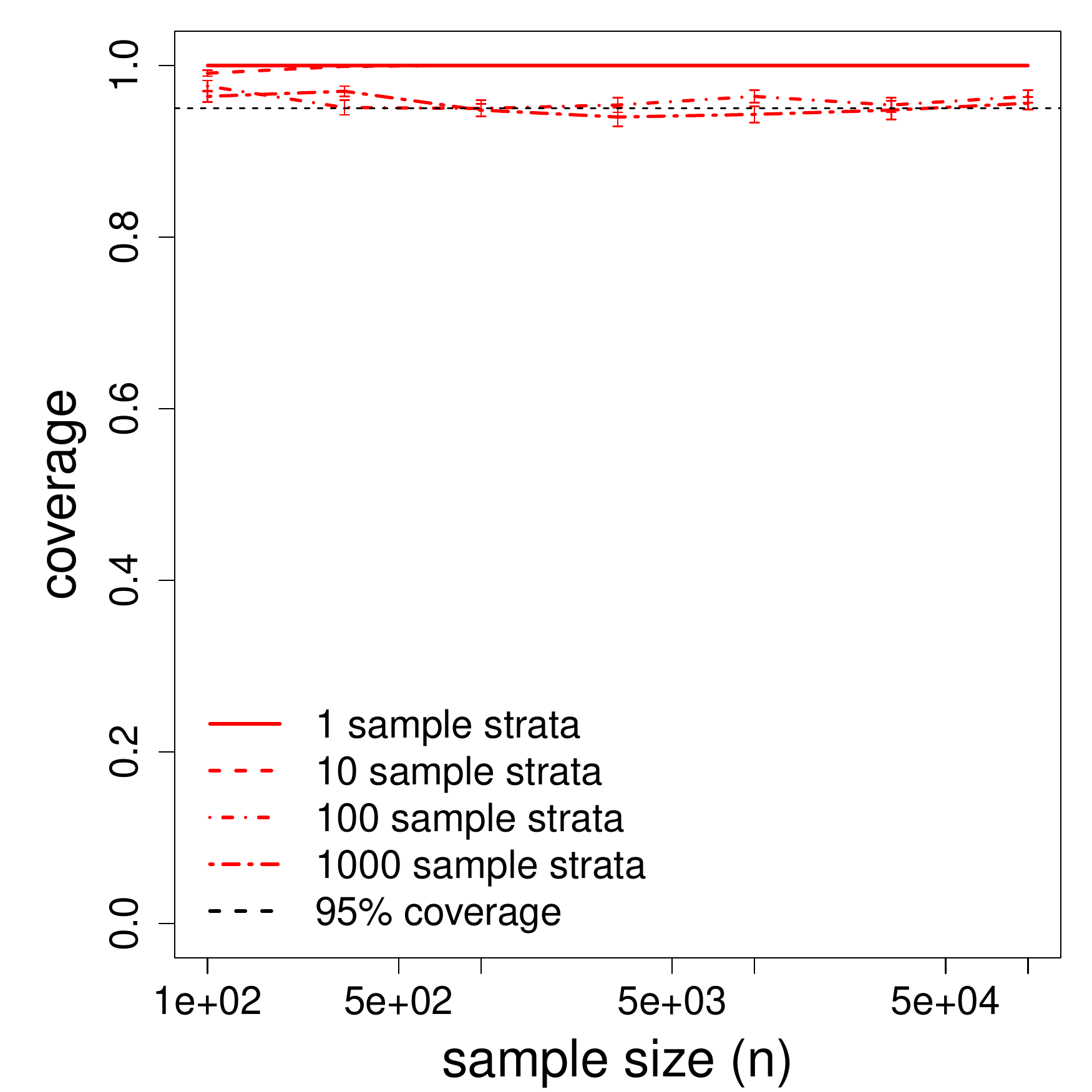}
\caption{Coverage versus number of samples per observation ($n_o$) for estimating 95\% CI using plugin method. We used AveAMSGrad with empirically optimal $C$ for learning rate defined as $\gamma_m=\frac{C}{\sqrt{m}}$), batch size $K=1$, $P=10$ features, strata size $S=20$, and probability of censoring $p_c=0.2$.}
\label{fig05}
\end{figure}

In Section~5.1, we gave a plugin approach for constructing a $(1-\alpha)\times100\%$ CI (see equations \eqref{u5} and \eqref{u6} in the main manuscript). Although there are ${n-1\choose s-1}$ possible  strata per each observation, in practice, we prefer to consider a random small fraction of them to estimate the standard error. Figure \ref{fig05} presents the coverage of $95\%$ CI constructed using the plugin approach with varying number of  strata sampled per observation (i.e., $n_o$). We considered 100 simulated datasets (based on the data generating mechanism explained in Section 6.1 with probability of censoring $p_c=0.2$) with 100 epochs to estimate $\tilde{\boldsymbol{\beta}}$. We observe sampling between 100 and 1000 strata per observation ($n_o$) gives near nominal coverage. This makes our computation tractable, especially for large datasets (e.g., for sample size $n=10000$, we only need 100 sample strata instead of going through all ${9999\choose 19}\approx 8.06\times 10^{58}$ sample strata).

\subsection{Practical data}
\label{secAG2}
We next compare our proposed framework BigSurvSGD with {\tt{coxph()}} using publicly available datasets. Note that the true parameter ($\boldsymbol{\beta}^*$) is unknown for these datasets. Therefore, we compare methods through: (1) estimation of hazards ratio with 95\% confidence interval and (2) concordance index. 

\subsubsection{Estimation of hazards ratio with 95\% confidence interval}
\label{secAG21}
Table \ref{Table1} presents estimates of hazard ratio with 95\% confidence interval for the Assay Of Serum Free Light Chain (FLCHAIN) dataset. This dataset is available in the  {\tt{survival}} package of R (Therneau, 2015). The aim of this study was to determine the association between serum free light chain and mortality \citep{Kyle2006}. We excluded variable "chapter" as it is the cause of death. We also removed observations with missing values. The resulting dataset has 6542 individuals with eight features: age, sex, the calendar year in which a blood sample was obtained (sample year), serum-free light chain, kappa portion (kappa), serum-free light chain, lambda portion (lambda), group of free light chain (FLC group), serum creatinine, and diagnosed with monoclonal gammopathy (0=no, 1=yes). We chose strata size $S=20$, mini-batch size $K=1$, $C=0.12$ in the learning rate defined as $\gamma_m=\frac{C}{\sqrt{m}}$, 100 epochs, and AMSGrad algorithm with average of estimate over iterates (AveAMSGrad). For the plugin method, we used $n_o=1000$ strata per observation to estimate the standard error. For the bootstrap method, we used $B=1000$ bootstrap resamples with 100 epochs. As we see, all three approaches give very close estimates of hazards ratio and 95\% confidence intervals for the all covariates' coefficients.

\begin{table}[!tb]
\centering
\caption{Estimates of hazards ratio and $95\%$ confidence interval for different covariates in dataset FLCHAIN ($n=6542$, $P=8$, $p_c=0.7$) with {\tt{coxph()}} and BigSurvSGD (AveAMSGrad). We use strata size $S=20$, batch size $K=1$, 100 epochs, and $C=0.12$ for the learning rate defined as $\gamma_m=\frac{C}{\sqrt{m}}$.}
\label{Table1}
\begin{tabular}{lccc}
 & {\tt{coxph()}} & BigSurvSGD-plugin & BigSurvSGD-bootstrap\\
 \hline
age & 1.107 (1.102, 1.113) & 1.107 (1.099, 1.114) & 1.107 (1.100, 1.113)\\
sex & 0.756 (0.688, 0.831) & 0.755 (0.674, 0.847) & 0.756 (0.673, 0.844)\\
sample year & 1.056 (1.017, 1.096) & 1.054 (1.015, 1.095) & 1.054 (1.008, 1.105)\\
kappa & 1.018 (0.952, 1.088) & 1.031 (0.939, 1.130) & 1.029 (0.891, 1.149)\\
lambda & 1.185 (1.122, 1.251) & 1.181 (1.111, 1.254) & 1.183 (1.097, 1.311)\\
FLC group & 1.057 (1.035, 1.079) & 1.056 (1.032, 1.080) & 1.056 (1.026, 1.090)\\
serum creatinine & 1.036 (0.941, 1.141) & 1.044 (0.932, 1.170) & 1.044 (0.854, 1.244)\\
mgus & 1.303 (0.792, 2.146) & 1.292 (0.817, 2.045) & 1.293 (0.857, 2.976)
\end{tabular}
\end{table}

\subsubsection{Concordance index}
\label{secAG22}
Table~\ref{Table3} presents the average concordance index with {\tt{coxph()}} and BigSurvSGD algorithms for different datasets. Details for the FLCHAIN dataset were given in the previous section. Dataset VETERAN is from a randomized trial of two treatment regimens for lung cancer \citep{Kalbfleisch2002}. This dataset is available in the {\tt{survival}} package of R \citep{Therneau2020}. This dataset has 137 individuals with six features: treatment (0=standard and 1=test), cell type (squamous, small, adeno, and large), Karnofsky performance score (100=good), age, prior therapy (0=no, 10=yes). Dataset German Breast Cancer Study Group (GBSG) is from a randomized 2 x 2 trial evaluating hormonal treatment and the duration of chemotherapy in node-positive breast cancer patients \citep{Schumacher1994}. This dataset is available in {\tt{mfp}} package in R \citep{Ambler2015}. This dataset has 686 individuals (women) with eight features: hormone therapy (0=no, 1=yes), age, menopausal status (0=no, 1=yes), tumor size, tumor grade ($1<2<3$), number of positive nodes, progesterone receptor, and estrogen receptor. Note that we used dummy variables for the categorical variables (e.g., tumor grade) with more than two categories. We chose strata size $S=20$, mini-batch size $K=1$, $C=0.12$ in the learning rate defined as $\gamma_m=\frac{C}{\sqrt{m}}$, 100 epochs, and AMSGrad algorithm with average of estimate over iterates (AveAMSGrad). Our proposed framework works at least as well as {\tt{coxph()}} for all datasets and all choices of strata size. The average concordance index increases by decreasing strata size (the highest concordance belongs to $S=2$). The explanation is straightforward: The concordance index is defined based on pairs of observations and is increased by choosing strata size as close as to $S=2$. Our framework also performs similarly to the one proposed by \citep{Kvamme2019} for dataset FLCHAIN with the same included individuals and features (average concordance index of 0.794 using our framework versus 0.793 using their framework).

\begin{table}[!tb]
\caption{Average concordance index with {\tt{coxph()}} and BigSurvSGD  (AveAMS-Grad) for three different datasets each with different sample size ($n$), number of features ($P$), and probability of censoring ($p_c$). We use strata size $S=20$, batch size $K=1$, 100 epochs, and $C=0.12$ for the learning rate defined as $\gamma_m=\frac{C}{\sqrt{m}}$.} 
\label{Table3}
\centering
\begin{tabular}{lcccc}
 & FLCHAIN & VETERAN & GBSG\\
& $n=6524$ & $n=137$ & $n=686$\\
& $P=8$ & $P=6$ & $P=8$\\
& $p_c=0.7$ & $p_c=0.07$ & $p_c=0.56$\\
\hline
{\tt{coxph()}} & 0.792 & 0.733 & 0.692\\
BigSurvSGD, S=2 & \textbf{0.794} & \textbf{0.745} & \textbf{0.696}\\
BigSurvSGD, S=5 & 0.793 & 0.744 & 0.695\\
BigSurvSGD, S=10 & 0.793 & 0.741 & 0.694\\
BigSurvSGD, S=20 & 0.793 & 0.739 & 0.693\\
BigSurvSGD, S=50 & 0.792 & 0.736 & 0.692
\end{tabular}
\end{table}

\bibliographystyle{biorefss}
\bibliography{Bibliography}

\begin{thebibliography}{99}

\bibitem[Aalen \emph{and others}(2008)Aalen, Borgan and
  Gjessing]{SurvPoint2008}
\textsc{Aalen, O., Borgan, O. and Gjessing, H.} (2008).
\newblock {\em Survival and Event History Analysis: A Process Point of View\/}.
  New York, USA: Springer.

\bibitem[Aggarwal(2018)Aggarwal]{NN}
\textsc{Aggarwal, C.~C.} (2018).
\newblock {\em Neural Networks and Deep Learning: A Textbook\/}. Switzerland:
  Springer.

\bibitem[Ambler \emph{and others}(2015)Ambler, Benner and Luecke]{Ambler2015}
\textsc{Ambler, G., Benner, A. and Luecke, S.} (2015).
\newblock {\em mfp: an R package for Multivariable Fractional Polynomials\/}.
\newblock R package version 1.5.2).

\bibitem[Andersen \emph{and others}(1993)Andersen, Borgan, Gill and
  Keiding]{Andersen1993}
\textsc{Andersen, P.~K., Borgan, O., Gill, R.~D. and Keiding, N.} (1993).
\newblock {\em Introduction to Empirical Processes and Semiparametric
  Inference\/}. NC, USA: Springer.

\bibitem[Austin and Steyerberg(2012)Austin and Steyerberg]{Peter2012}
\textsc{Austin, P.~C. and Steyerberg, E.~W.} (2012).
\newblock Interpreting the concordance statistic of a logistic regression
  model: relation to the variance and odds ratio of a continuous explanatory
  variable.
\newblock {\em BMC Medical Research Methodology\/}~\textbf{12}(82), 1--8.

\bibitem[Bach and Moulines(2011)Bach and Moulines]{Eric2011}
\textsc{Bach, F.~R. and Moulines, E.} (2011).
\newblock Non-asymptotic analysis of stochastic approximation algorithms for
  machine learning.
\newblock {\em Neural Information Processing Systems (NeurIPS)\/}.

\bibitem[Bender \emph{and others}(2005)Bender, Augustin and
  Blettner]{Bender2005}
\textsc{Bender, R., Augustin, T. and Blettner, M.} (2005).
\newblock Generating survival times to simulate cox proportional hazards
  models.
\newblock {\em Statistics in medicine\/}~\textbf{24}(11), 1713--1723.

\bibitem[Bickel and Freedman(1981)Bickel and Freedman]{Bickel1981}
\textsc{Bickel, P.~J. and Freedman, D.~A.} (1981).
\newblock Some asymptotic theory for the bootstrap.
\newblock {\em The Annals of Statistics\/}~\textbf{9}(6), 1196--1217.

\bibitem[Blom(1976)Blom]{Gunnar1976}
\textsc{Blom, G.} (1976).
\newblock Some properties of incomplete u-statistics.
\newblock {\em Biometrika\/}~\textbf{63}(3), 573--580.

\bibitem[Bottou(2010)Bottou]{Bottou2010}
\textsc{Bottou, L.} (2010).
\newblock Large-scale machine learning with stochastic gradient descent.
\newblock In: Lechevallier, Yves and Saporta, Gilbert (editors), {\em
  Proceedings of COMPSTAT'2010\/}. Paris, France: Springer. pp.\  177--186.

\bibitem[Boyd \emph{and others}(2010)Boyd, Parikh, Chu, Peleato and
  Eckstein]{Boyd2010}
\textsc{Boyd, S., Parikh, N., Chu, E., Peleato, B. and Eckstein, J.} (2010).
\newblock Distributed optimization and statistical learning via the alternating
  direction method of multipliers.
\newblock {\em Foundations and Trends in Machine Learning\/}~\textbf{3}(1),
  1--122.

\bibitem[Boyd and Vandenberghe(2004)Boyd and Vandenberghe]{Boyd2004}
\textsc{Boyd, S. and Vandenberghe, L.} (2004).
\newblock {\em Convex Optimization\/}. Cambridge, UK: Cambridge University
  Press.

\bibitem[Breslow and Day(1980)Breslow and Day]{Breslow1980}
\textsc{Breslow, N.~E. and Day, N.~E.} (1980).
\newblock Statistical methods in cancer research.vol 1. the analysis of
  case-control studies.
\newblock {\em International Agency for Research on Cancer\/}~\textbf{1}(32),
  5--338.

\bibitem[Chapfuwa \emph{and others}(2018)Chapfuwa, Tao, Li, Page, Goldstein,
  Carin and Henao]{Chapfuwa2018}
\textsc{Chapfuwa, P., Tao, C., Li, C., Page, C., Goldstein, B., Carin, L. and
  Henao, R.} (2018).
\newblock Adversarial time-to-event modeling.
\newblock {\em Proceedings of 35th International Conference on Machine
  Learning\/}.

\bibitem[Ching \emph{and others}(2018)Ching, Zhu and Garmire]{Ching2018}
\textsc{Ching, T., Zhu, X. and Garmire, L.~X.} (2018).
\newblock Cox-nnet: An artificial neural network method for prognosis
  prediction of high-throughput omics data.
\newblock {\em Plus Computational Biology\/}~\textbf{14}(4), 1--18.

\bibitem[Cox(1972)Cox]{Cox1972}
\textsc{Cox, D.~R.} (1972).
\newblock Regression models and life-tables.
\newblock {\em Journal of the Royal Statistical Society. Series B
  (Methodological)\/}~\textbf{34}(2), 187--220.

\bibitem[Duchi \emph{and others}(2011)Duchi, Hazan and Singer]{Duchi2011}
\textsc{Duchi, J., Hazan, E. and Singer, Y.} (2011).
\newblock Adaptive subgradient methods for online learning and stochastic
  optimization.
\newblock {\em Journal of Machine Learning Researc\/}~\textbf{12}(7),
  2121--2159.

\bibitem[Efron(1979)Efron]{Efron1979}
\textsc{Efron, B.} (1979).
\newblock Bootstrap methods: another look at the jackknife.
\newblock {\em The Annals of Statistics\/}~\textbf{7}(1), 1--26.

\bibitem[Fang \emph{and others}(2018)Fang, Xu and Yang]{Fang2018}
\textsc{Fang, Y., Xu, J. and Yang, L.} (2018).
\newblock Online bootstrap confidence intervals for the stochasticgradient
  descent estimator.
\newblock {\em Journal of Machine Learning Research\/}~\textbf{1}(19),
  3053--3073.

\bibitem[Gaber \emph{and others}(2005)Gaber, Zaslavsky and
  Krishnaswamy]{Gaber2005}
\textsc{Gaber, M.~M., Zaslavsky, A. and Krishnaswamy, S.} (2005).
\newblock Mining data streams: A review.
\newblock {\em ACM Digital Library\/}~\textbf{34}(2), 18--26.

\bibitem[Honore and Powell(1994)Honore and Powell]{Honore1994}
\textsc{Honore, B.~E. and Powell, J.~L.} (1994).
\newblock Pairwise difference estimators of censored and truncated regression
  models.
\newblock {\em Journal of Econometrics\/}~\textbf{64}(1), 241--278.

\bibitem[Kalbfleisch and Prentice(2002)Kalbfleisch and
  Prentice]{Kalbfleisch2002}
\textsc{Kalbfleisch, J.~D. and Prentice, R.~L.} (2002).
\newblock {\em The Statistical Analysis of Failure Time Data, Second
  Edition\/}. NJ, USA: John Wiley and Sons.

\bibitem[Katzman \emph{and others}(2017)Katzman, Shaham, Cloninger, Bates,
  Jiang and Kluger]{Katzman2017}
\textsc{Katzman, J.~L., Shaham, U., Cloninger, A., Bates, J., Jiang, T. and
  Kluger, Y.} (2017).
\newblock Deepsurv: Personalized treatment recommender system using a cox
  proportional hazards deep neural network.
\newblock eprint arXiv:1606.00931v3.

\bibitem[Kingma and Ba(2014)Kingma and Ba]{Kingma2014}
\textsc{Kingma, D.~P. and Ba, J.~L.} (2014).
\newblock Adam: A method for stochastic optimization.
\newblock eprint arXiv:1412.6980.

\bibitem[Klein and Moeschberger(2003)Klein and Moeschberger]{Klein2003}
\textsc{Klein, J.~P. and Moeschberger, M.~L.} (2003).
\newblock {\em Survival Analysis: Techniques for Censored and Truncated
  Data\/}. New York, USA: Springer.

\bibitem[Koch \emph{and others}(2015)Koch, Zemel and
  Salakhutdinov]{Gregory2015}
\textsc{Koch, G., Zemel, R. and Salakhutdinov, R.} (2015).
\newblock Siamese neural networks for one-shot image recognition.
\newblock {\em Proceedings of the 32nd International Conference on Machine
  Learning\/}~\textbf{37}.

\bibitem[Kvamme \emph{and others}(2019)Kvamme, Borgan and Scheel]{Kvamme2019}
\textsc{Kvamme, H., Borgan, O. and Scheel, I.} (2019).
\newblock Time-to-event prediction with neural networks and cox regression.
\newblock {\em Journal of Machine Learning Research\/}~\textbf{20}(129), 1--30.

\bibitem[Kyle \emph{and others}(2006)Kyle, Therneau, Rajkumar, Larson, Plevak,
  Offord, Dispenzieri, Katzmann and Melton]{Kyle2006}
\textsc{Kyle, R., Therneau, T., Rajkumar, S.~V., Larson, D., Plevak, M.,
  Offord, J., Dispenzieri, A., Katzmann, J. and Melton, L.~J.} (2006).
\newblock Prevalence of monoclonal gammopathy of undetermined significance.
\newblock {\em New England Journal of Medicine\/}~\textbf{354}, 1362--1369.

\bibitem[Lee and Go(1997)Lee and Go]{Lee1997}
\textsc{Lee, E.~T. and Go, O.~T.} (1997).
\newblock Survival analysis in public health research.
\newblock {\em Annual Reviews in Public Health\/}~\textbf{18}, 105--134.

\bibitem[Mittal and Madigan(2014)Mittal and Madigan]{MITTAL2014}
\textsc{Mittal, S. and Madigan, D.} (2014).
\newblock High-dimensional, massive sample-size cox proportional hazards
  regression for survival analysis.
\newblock {\em Biostatistics\/}~\textbf{15}(2), 207--221.

\bibitem[Polyak and Juditsky(1992)Polyak and Juditsky]{Polyak1992}
\textsc{Polyak, B.~T. and Juditsky, A.~B.} (1992).
\newblock Acceleration of stochastic approximation by averaging.
\newblock {\em SIAM Journal on Control and Optimization\/}~\textbf{30}(4),
  838--855.

\bibitem[Raghupathi and Raghupathi(2014)Raghupathi and
  Raghupathi]{Raghupathi2014}
\textsc{Raghupathi, W. and Raghupathi, V.} (2014).
\newblock Big data analytics in healthcare: promise and potential.
\newblock {\em Journal of Health Information Science and
  Systems\/}~\textbf{2}(3), 1--10.

\bibitem[Raykar \emph{and others}(2008)Raykar, Steck, Krishnapuram,
  Dehing-Oberije and Lambin]{Raykar2008}
\textsc{Raykar, V.~C., Steck, H., Krishnapuram, B., Dehing-Oberije, C. and
  Lambin, P.} (2008).
\newblock On ranking in survival analysis: Bounds on the concordance index.
\newblock {\em Neural Information Processing Systems (NeurIPS)\/}.

\bibitem[Reddi \emph{and others}(2019)Reddi, Kale and Kumar]{Sashank2019}
\textsc{Reddi, S.~J., Kale, S. and Kumar, S.} (2019).
\newblock On the convergence of adam and beyond.
\newblock eprint arXiv:1904.09237.

\bibitem[Ruder(2016)Ruder]{Ruder2016}
\textsc{Ruder, S.} (2016).
\newblock An overview of gradient descent optimization algorithms.
\newblock eprint arXiv:1609.04747.

\bibitem[Ruppert(1988)Ruppert]{Ruppert1988}
\textsc{Ruppert, D.} (1988).
\newblock Efficient estimations from a slowly convergent robbins-monro process.
\newblock {\em Technical Report}, School of Operations Research and Industrial
  Engineering, Cornell University, New York, USA.

\bibitem[Schober and Vetter(2018)Schober and Vetter]{Schober2018}
\textsc{Schober, P. and Vetter, T.~R.} (2018).
\newblock Survival analysis and interpretation of time-to-event data: The
  tortoise and the hare.
\newblock {\em Anesthesia \& Analgesia\/}~\textbf{127}(3), 792--798.

\bibitem[Schumacher \emph{and others}(1994)Schumacher, Bastert, H, Hubner,
  Olschewski, Sauerbrei, Schmoor, Beyerle, Neumann and
  Rauschecker]{Schumacher1994}
\textsc{Schumacher, M., Bastert, G., H, H.~Bojar, Hubner, K., Olschewski, M.,
  Sauerbrei, W., Schmoor, C., Beyerle, C., Neumann, R.~L. and Rauschecker,
  H.~F.} (1994).
\newblock Randomized 2 x 2 trial evaluating hormonal treatment and the duration
  of chemotherapy in node-positive breast cancer patients. german breast cancer
  study group.
\newblock {\em Journal of Clinical Oncology\/}~\textbf{12}(10), 2086--2093.

\bibitem[Sherman(1994)Sherman]{Sherman1994}
\textsc{Sherman, R.~P.} (1994).
\newblock Maximal inequalities for degenerate u-processes with applications to
  optimization estimators.
\newblock {\em The Annals of Statistics\/}~\textbf{22}(1), 439--459.

\bibitem[Simon \emph{and others}(2011)Simon, Friedman, Hastie and
  Tibshirani]{Simon2011}
\textsc{Simon, N., Friedman, J., Hastie, T. and Tibshirani, R.} (2011).
\newblock Regularization paths for cox’s proportional hazards model via
  coordinate descent.
\newblock {\em Journal of Statistical Software\/}~\textbf{39}(5), 1--13.

\bibitem[Tarkhan and Simon(2020)Tarkhan and Simon]{Tarkhan2020}
\textsc{Tarkhan, A. and Simon, N.} (2020).
\newblock Big survival data analysis via stochastic gradient descent.
\newblock \url{https://github.com/atarkhan/bigSurvSGD}.

\bibitem[Therneau and Lumley(2019)Therneau and Lumley]{coxph2019}
\textsc{Therneau, T.~M. and Lumley, T.} (2019).
\newblock coxph: Core survival analysis routines.
\newblock R package version 2.44-1.1.

\bibitem[Therneau \emph{and others}(2020)Therneau, Lumley, Elizabeth and
  Cynthia]{Therneau2020}
\textsc{Therneau, T.~M., Lumley, T., Elizabeth, A. and Cynthia, C.} (2020).
\newblock {\em survival: An R package for Survival Analysis\/}.
\newblock R package version 3.1-12).

\bibitem[Toulis and Airoldi(2017)Toulis and Airoldi]{Toulis2017}
\textsc{Toulis, P. and Airoldi, E.~M.} (2017).
\newblock Asymptotic and finite-sample properties of estimators based on
  stochastic gradients.
\newblock {\em The Annals of Statistics\/}~\textbf{45}(4), 1694--1727.

\bibitem[van~der Vaart(2000)van~der Vaart]{Van2000}
\textsc{van~der Vaart, A.~W.} (2000).
\newblock {\em Asymptotic Statistics\/}. Cambridge, UK: Cambridge University
  Press.

\bibitem[Wei(1992)Wei]{Wei1992}
\textsc{Wei, L.~J.} (1992).
\newblock The accelerated failure time model: a useful alternative to the cox
  regression model in survival analysis.
\newblock {\em Statistics in Medicine\/}~\textbf{11}, 1871--1879.

\bibitem[Zethelius \emph{and others}(2008)Zethelius, Berglund, Sundstrom, ,
  Ingelsson, Basu, Larsson, Venge,  and Arnlov]{Bjorn2008}
\textsc{Zethelius, B., Berglund, L., Sundstrom, J., , Ingelsson, E., Basu, S.,
  Larsson, A., Venge, P.,  and Arnlov, J.} (2008).
\newblock Use of multiple biomarkers to improve the prediction of death from
  cardiovascular causes.
\newblock {\em The New England Journal of Medicine\/}~\textbf{358}(20),
  2107--2116.

\end{thebibliography}

\end{document}